%% file: tac_arxiv.tex
\documentclass[12pt]{article}
\usepackage[margin=2.5cm]{geometry}

\input{./Texts/Initialize}
\begin{document}
\title{\bf Analysis and Design of Vehicle Platooning Operations on Mixed-Traffic Highways}
\author{Li~Jin, Mladen~\v Ci\v ci\'c, Karl H. Johansson, and Saurabh~Amin
\thanks{
This work was supported by NYU Tandon School of Engineering, C2SMART University Transportation Center, US National Science Foundation CAREER Award CNS-1453126, US Air Force Office of Scientific Research, EU's Horizon 2020 Marie Sk{\l}odowska-Curie 674875, FFI VINNOVA 2014-06200, Swedish Research Council, Swedish Foundation for Strategic Research, and Knut and Alice Wallenberg Foundation. Xi Xiong and Teze Wang at NYU Tandon School of Engineering contributed to the simulation. The authors greatly appreicate the feedback from the anonymous reviewers and from the associate editor.}
\thanks{L. Jin is with the Tandon School of Engineering, New York University, Brooklyn, New York, USA.
M. \v Ci\v ci\'c and K. H. Johansson are with the School of Electrical Engineering and Computer Science, KTH Royal Institute of Technology, Stockholm, Sweden.
S. Amin is with the Laboratory for Information and Decision Systems and the Department of Civil and Environmental Engineering, Massachusetts Institute of Technology, Cambridge, Massachusetts, USA. (e-mails: lijin@nyu.edu, cicic@kth.se, kallej@kth.se, amins@mit.edu).}
}
\maketitle

\input{./Texts/0_Abstract}
\input{./Texts/1_Introduction}
\input{./Texts/2_Modeling}
\input{./Texts/3_Analysis}

\input{./Texts/4_Control}
\input{./Texts/4.5_Simulate}
\input{./Texts/5_Conclusion}

\bibliographystyle{IEEEtran}
\bibliography{Bibliography_LJ,Bibliography_KHJ}
\end{document}

%% file: Texts/Initialize.tex
\usepackage{subfigure}
\usepackage{amsmath}

\usepackage{amsthm} 
\allowdisplaybreaks
\usepackage[usenames, dvipsnames]{color}
\usepackage{footnote}
\usepackage{algorithm}
\usepackage{graphicx}
\usepackage{amssymb}
\usepackage{amsbsy}
\usepackage{array}
\usepackage{longtable}
\usepackage{epstopdf}
\usepackage{pbox}
\usepackage{url}
\usepackage{breqn}
\usepackage{mathrsfs}
\usepackage{multicol}
\usepackage{supertabular}
\usepackage{enumerate}
\usepackage{hyperref}
\usepackage{cite}
\usepackage{mathtools}
\DeclarePairedDelimiter{\ceil}{\lceil}{\rceil}
\usepackage{booktabs}
\newtheorem{dfn}{Definition}
\newtheorem{thm}{Theorem}

\newtheorem{rem}{Remark}
\newtheorem{lmm}{Lemma}
\newtheorem{asm}{Assumption}

%% file: Texts/0_Abstract.tex
\begin{abstract}
Platooning of connected and autonomous vehicles (CAVs) has a significant potential for throughput improvement. However, the interaction between CAVs and non-CAVs may limit the practically attainable improvement due to platooning. To better understand and address this limitation, we introduce a new fluid model of mixed-autonomy traffic flow and use this model to analyze and design platoon coordination strategies. We propose tandem-link fluid model that considers randomly arriving platoons sharing highway capacity with non-CAVs. We derive verifiable conditions for stability of the fluid model by analyzing an underlying M/D/1 queuing process and establishing a Foster-Lyapunov drift condition for the fluid model. These stability conditions enable a quantitative analysis of highway throughput under various scenarios. The model is useful for designing platoon coordination strategies that maximize throughput and minimize delay. Such coordination strategies are provably optimal in the fluid model and are practically relevant. We also validate our results using standard macroscopic (cell transmission model, CTM) and microscopic (Simulation for Urban Mobility, SUMO) simulation models.
\end{abstract}


%% file: Texts/1_Introduction.tex
\section{Introduction}
\label{sec_intro}

\begin{figure}[hbt]
\centering
\subfigure[Without inter-platoon coordination, the exiting traffic can be blocked by local congestion induced by platoons.]{
\centering
\includegraphics[width=0.4\textwidth]{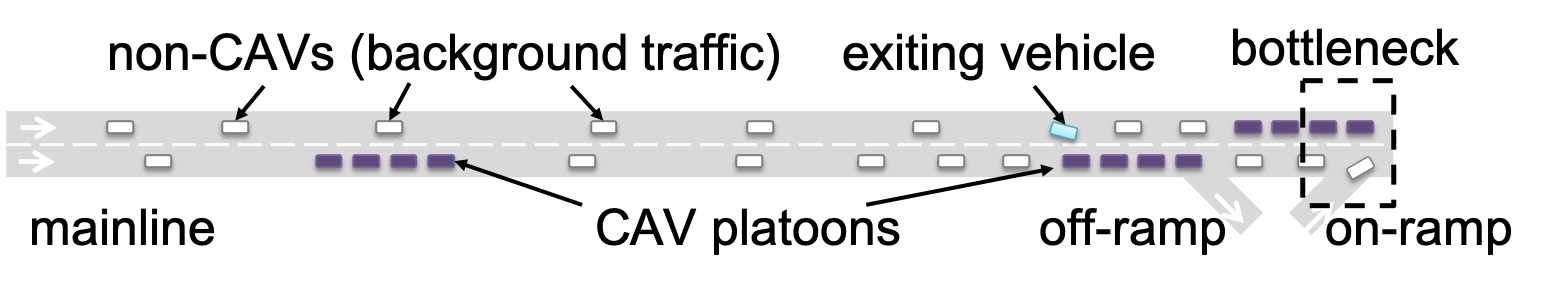}
\label{cartoon_nocontrol}
}
\subfigure[Inter-platoon coordination can help mitigate the ``moving-obstacle'' effect and maintain free flow.]{
\centering
\includegraphics[width=.4\textwidth]{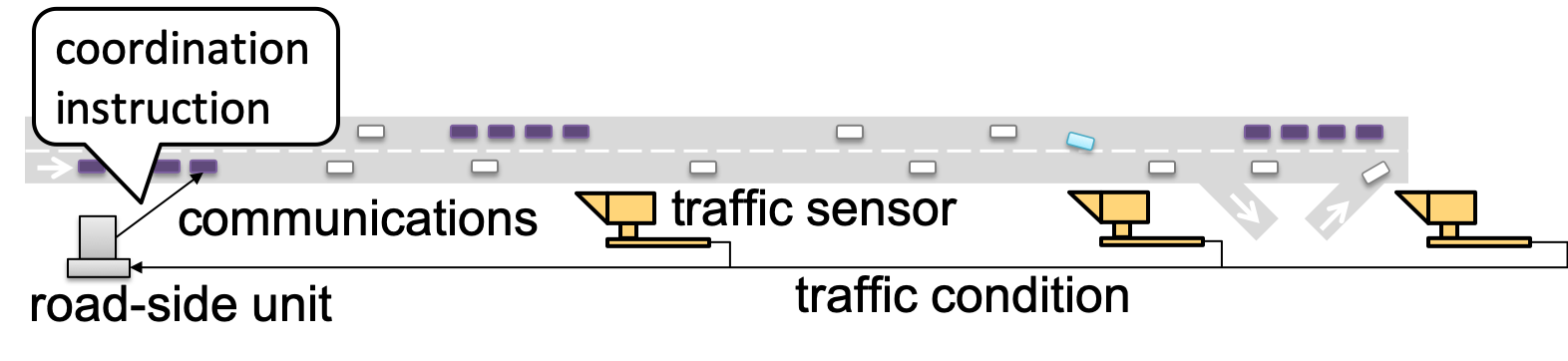}
\label{cartoon_control}
}
\caption{Traffic scenarios without (top) and with (bottom) inter-platoon coordination.}
\end{figure}

Platooning of connected and autonomous vehicles (CAVs) has the potential for significant throughput improvement~\cite{hor+var00,litman2017autonomous}. The idea of automatically regulating strings of vehicles is well-known~\cite{lev+ath66,swaroop1996string}, and several experimental studies in real-world traffic conditions have been conducted in the past decades~\cite{chang1993automated,naus+10,bess+16procieee,tsugawa16}. 
Important progress has been made in vehicle-level control design~\cite{naus+10,coogan15interconnected,bess+16procieee} and system-level simulations \cite{talebpour16}.
However, we still lack both link- or network-level models for evaluating the impact of platooning on highway traffic and coordinating the movement of platoons on mixed-traffic highways.
In particular, although platooned CAVs have smaller inter-vehicle spacing, uncoordinated and randomly arriving CAV platoons may act as ``moving obstacles'' and result in recurrent local congestion, especially at bottlenecks; see Fig.~\ref{cartoon_nocontrol} for illustration. 
Appropriate inter-platoon coordination, such as regulating the headways between platoons and managing the platoon sizes, is essential for realizing the full benefits of platooning in mixed-traffic conditions.
%

In this paper, we consider a platoon coordination problem over a generic highway section as shown in Fig.~\ref{cartoon_nocontrol}.
The highway section has a downstream on-ramp, which forms a geometric bottleneck.
If more traffic is arriving at the bottleneck than the bottleneck can discharge, then a queue would be growing there.
The bottleneck may impede free flow of traffic, and consequently uncoordinated platoons may start inducing congestion (via local queuing), thereby impacting the upstream off-ramp.
One of the objectives of platoon coordination is ensure that the potential throughput gain due to platooning is not limited by such local congestion.
Inter-platoon coordination can be achieved by (i) regulating the times at which platoons arrive at the bottleneck by specifying the average or reference speed of platoons over the highway section (called \emph{headway regulation}) and (ii) maintaining desirable platoon sizes through split or merge maneuvers (called \emph{size management}).
Implementation of such coordination strategies can be facilitated by traffic sensors that collect real-time traffic information and road-side units that send coordination instructions to each platoon (Fig.~\ref{cartoon_control}); such capabilities are already available in modern transportation systems \cite{duret2019hierarchical}.

To model the interaction between platoons and background traffic, we present a multiclass tandem-link fluid model.
Fluid models are standard for highway bottleneck analysis \cite[Ch. 2]{newell13} and allow tractable analysis \cite{kulkarni97}.
Our model specifically captures two important features of platooning operations.
First,  platoons travel in a clustered manner (i.e. with small between-vehicle spacings and relatively large inter-platoon headways), whereas non-CAVs do not follow such a configuration. 
In our model, random platoon arrivals are modeled as Poisson jumps in the traffic queue at the bottleneck.
This model captures the inherent randomness in platooning operations \cite{bess+16procieee}. Second, platoons share highway capacity with the background traffic and may act as ``moving obstacles''; our model captures this interaction by constraining the sum of discharge rates of platoons and background traffic at the bottleneck with an overall capacity.


For comparison, Table~\ref{tab_existing} lists various models that researchers have studied for traffic flow with CAV platoons under various contexts.
\begin{table}[hbt]
\caption{Models for traffic with platooning.}
\label{tab_existing}
\centering
\begin{tabular}{@{}cclc@{}}
\toprule
Scale                  & Model & Control actions & References      \\ \midrule
Vehicle & 
\begin{tabular}{@{}l@{}}Vehicle dynamic,\\ car following\end{tabular}  
& \begin{tabular}{@{}l@{}}Throttle, break,\\ steering\end{tabular} 
&\begin{tabular}{@{}c@{}} \cite{li2017dynamical,smith2019balancing}	\\
\cite{zhou2019stabilizing,duret2019hierarchical}
\end{tabular}    \\ \hline
Intersection & Discrete queuing  & 
 \begin{tabular}{@{}l@{}}Signal timing,\\ CAV coordination\end{tabular}
 &\cite{lioris17,miculescu2019polling}          \\ \hline
Segment              &  
\begin{tabular}{@{}l@{}}Partial differential\\
equation\end{tabular} 
& 
\begin{tabular}{@{}l@{}}Speed regulation,\\ platoon merge/split, \\  lane management\end{tabular} & \cite{keim+17,vcivcic2019coordinating}      \\ \hline
Link            & Fluid  & 
\begin{tabular}{@{}l@{}}Headway regulation,\\ platoon merge/split, \\  \end{tabular} & this paper 
 \\ \hline
Network        &   Static     & Routing  & \cite{lazar2018routing,mehr2018can}      \\ \bottomrule
\end{tabular}
\end{table}
%
%
Fluid models are particularly suitable for link-level analysis and control design due to the following advantages.
Compared with discrete queuing models, fluid models do not do not track individual vehicles; instead, only the aggregate flow is required. This considerably simplifies the state space and the system dynamics. Compared with PDE models, fluid models entail smaller computational requirements and enable tractable analysis. Compared with static models, fluid models allow real-time control design rather than long-term decisions such as day-to-day traffic assignment.

Two major differences are also worth noting here.
First, the behavior of discrete queuing models (e.g. M/M/1) is closely related to their fluid counterparts \cite{dai1995stability}. 
As opposed to fully discrete queuing models, our model considers platoon as discrete Poisson arrivals and non-CAVs as a continuous inflow.
This is motivated by practical situations where the arrival rate of CAV platoons (less than five per minute) is much lower than that of non-CAVs (50--100 per minute).
Second, fluid models share some common features with PDE models \cite{yu2019traffic} or their discretization (e.g. cell transmission model \cite{cicic2018traffic}). 
In contrast to PDE models, the fluid model retains the queuing delay due to demand-capacity imbalance but does not capture the evolution of congestion waves.
Still, both the traffic fluid model and classical traffic flow models are based on conservation laws. Some authors showed that these two types models lead to equivalent results in traffic network optimization \cite{qian2012system}. This paper also demonstrates the consistency between the fluid model and more detailed models via simulation (Section~\ref{sub_simulation}).

Using the fluid model, we study the throughput of the highway section with uncoordinated platoons\footnote{In this paper, ``coordination'' refers to the link-level coordination of multiple platoons. We acknowledge that ``coordination'' is also used to refer to the CACC applied to individual vehicles within one platoon \cite{vandehoef2018fuel}.}.
We utilize an M/D/1 queuing characterization of the fluid model and establish a Foster-Lyapunov drift condition for stability of the fluid model.
This leads to an easy-to-check sufficient condition for bounded queues of the uncontrolled system (Theorem~\ref{thm_sufficient}) which relies on the general stability/convergence theory of Markov processes \cite{meyn1993stability}.
We also derive explicit lower and upper bounds for throughput in the uncoordinated scenario (Theorem~\ref{thm_bounds}).
These results also contribute to the literature on fluid queuing systems \cite{mitra1988stochastic,kulkarni97,cassandras02,kroese2001joint,jin2018stability}


We also design a class of platoon coordination strategies that realize the full potential of platooning for throughput improvement.
The control actions that we consider include regulating inter-platoon headways and splitting platoons into shorter platoons.
In terms of the fluid model, these regulation strategies are formulated as control laws regulating the arrival process of platoons with the knowledge of the arrival times of previous platoons. In practice, such knowledge can be obtained by tracking the movement of existing platoons on the highway section via vehicle-to-infrastructure communications \cite{vandehoef2018fuel}.
We prove that a set of control laws (Theorem~\ref{thm_optimal}) stabilize the system in a fairly strong sense (bounded moment-generating function and exponential convergence to steady-state distribution) if and only if the total inflow is less than the capacity.
Thus, they are optimal in the sense of throughput maximization and delay minimization.
Intuitively, these control laws coordinate the movement of platoons so that they arrive at the bottleneck with relatively evenly distributed headways;
thus, queuing delay is absorbed en-route, congestion at the bottleneck does not build up, and spillback is eliminated.


Note that the control actions considered in this paper are related to but different from a class of longitudinal/lateral Cooperative Adaptive Cruise Control (CACC) capabilities. We focus on link-level decision variables including (i) the reference speed or the average speed that a traffic manager would recommend a platoon to take over a highway section and (ii) the decision whether to maintain or split a platoon over a highway section, both of which are concerned with a typical scale of 10 km or 10 minutes. In the context of vehicle-level CACC, however, the real-time speed is dynamically regulated to maintain string stability \cite{li2017dynamical} or to mitigate local stop-and-go behavior \cite{stern2018dissipation,wu2018stabilizing,zhou2019stabilizing}, which involve much finer space and time resolutions. The objective of CACC design is to optimize microscopic driving behavior and regulate congestion waves, while our control objective is to reduce the congestion due to randomness in platoon arrivals and develops coordination strategies that work effectively in the presence of such randomness.

Finally, we discuss how the fluid model-based results can be translated to implementable actions for actual CAV platoons and validate the proposed coordination strategies in standard simulation environments.
We implement the proposed strategies in a macroscopic (cell transmission model, CTM \cite{daganzo94}) and a microscopic (Simulation for Urban Mobility, SUMO \cite{krajzewicz2002sumo}) simulation model.
Results indicate that the proposed strategy effectively and consistently improves travel times in both simulation models.
In spite of multiple simplifications that our modeling approach makes, simulation results suggest that the theoretically optimal headway regulation strategy attains more than 80\% of the improvement attained by the simulation-optimal strategies.

The main contributions of this paper include:
\begin{enumerate}[(i)]
\item A novel fluid model a highway with randomly arriving platoons and constant inflow of background traffic that captures the queuing and throughput loss due to interaction between various traffic classes;
\item A set of easily checkable conditions for stability of the fluid model derived by combining ideas from queuing theory (mainly M/D/1 process) and the theory of stability of continuous-time, continuous-state Markov processes. These conditions enable quantitative analysis of throughput of mixed-autonomy highways;
\item A set of platooning strategies that regulate movement of platoons to attain maximum throughput as well as minimum delay;
\item Validation of the fluid model-based analysis and design results via simulation of macroscopic and microscopic models.
\end{enumerate}

The rest of this article is organized as follows.
Section~\ref{sec_modeling} introduces the fluid model.
Section~\ref{sec_analysis} presents throughput analysis based on stability conditions.
Section~\ref{sec_control} discusses a class of optimal control strategies.
Section~\ref{sec_ctm} presents the implementation and validation in simulation environments.
Section~\ref{sec_conclude} summarizes the main results and mentions several directions for future work.

%% file: Texts/2_Modeling.tex
\section{Modeling and problem definition}
\label{sec_modeling}

In this section, we introduce our stochastic fluid model for highway bottlenecks with mixed traffic consisting of both CAVs and non-CAVs (Section~\ref{sub_modeling}) and formally define the problems that we study in the rest of this paper (Section~\ref{sub_preliminaries}).

\subsection{Stochastic fluid queuing model}
\label{sub_modeling}

\begin{figure}[hbt]
\centering
\subfigure[2-link highway section and main parameters.]{
\centering
\includegraphics[width=0.4\textwidth]{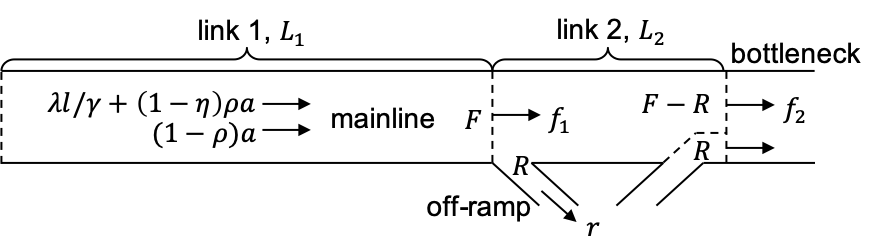}
\label{capacities}
}
\subfigure[Fluid model for system in Fig.~\ref{capacities}.]{
\centering
\includegraphics[width=0.3\textwidth]{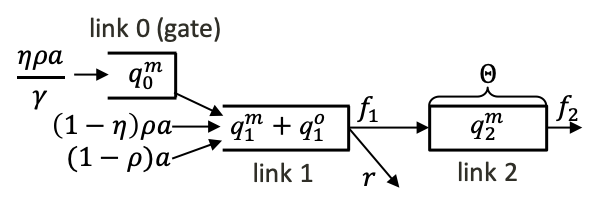}
\label{fig_tandem}
}
\caption{Fluid queuing system for highway section with CAVs and non-CAVs.}
\end{figure}

Consider a highway section with a downstream bottleneck and an off-ramp, as illustrated in Fig.~\ref{capacities}.  We model the highway as a tandem-link fluid queuing system in Fig.~\ref{fig_tandem}.
Link 1 has a \emph{mainline capacity} $F$ (veh/hr) and an \emph{off-ramp capacity} $R$ (veh/hr). Link 2 has an on-ramp with capacity $R$, which creates a bottleneck with capacity $F-R$.\footnote{For simplicity, we assume the on-ramp capacity to be identical to the off-ramp capacity.}
To model the effect of coordinated platooning operations (see Fig.~\ref{cartoon_control}), we introduce a ``virtual link'' 0 upstream to link 1, which we refer as the \emph{gate}. In our model, this gate can temporarily hold platoons and control their rate of release in order to regulate the downstream traffic flow. Furthermore, the storage space in link 2 (i.e., the maximum queue length that it can admit before the traffic spills over to link 1) is finite, modeled as a \emph{buffer} with space $\Theta$ (veh). On the other hand, links 0 and 1 are assumed to have an infinite buffer space. 
However, the model can be extended to other typical road configurations \cite{qian2012system}, which allows extension of our approach to more general settings. Note that our model does not account for details such as car following and lane changing.

The highway section is subject to a total demand $a>0$ (veh/hr). This demand comprises of $\rho a$ amount of mainline  demand that is discharged to link 2, and $(1-\rho)a$ amount of demand exiting through the off-ramp. Out of the $\rho a$ mainline demand, CAVs traveling in platoons amount to a fraction $\eta\in[0,1]$, and the remaining $(1-\eta)$ fraction is comprised of non-CAV traffic. In this demand pattern, CAVs only contribute to the mainline traffic (with demand $\eta\rho a$), and the off-ramp demand $(1-\rho)a$ is entirely comprised of non-CAVs.\footnote{Our setup can be extended to a more general case of CAV platoons that are bound to different destinations.} We call $\rho$ the \emph{mainline ratio} and $\eta$ the \emph{platooning ratio}. 
We refer to the total non-CAV demand $(1-\eta\rho)a$ as \emph{background traffic}.
For CAV platoons, the mean inter-arrival time is much greater than minimal inter-arrival times, so we consider platoons as discrete arrivals. For non-CAVs, however, the mean and the minimal inter-arrival times are in the same order, so we consider background traffic as continuous fluid.
We consider a Poisson process rather than a general renewal process, since Poisson processes are standard models for random arrivals in transportation systems \cite{newell13} and ensure that the model is Markovian \cite{benaim15}.

CAV platoons arrive according to a Poisson process of rate $\lambda$ (platoons per hr), which is given by
\begin{align}
    \lambda:=\frac{\eta\rho a}{l},
    \label{eq_lambda}
\end{align}
where $l$ is the number of CAVs in each platoon. Typical values of $l$ are between 2 and 10 \cite{bess+16procieee}.
For ease of presentation, we consider homogeneous platoon lengths. Non-homogeneous platoon lengths can be modeled as jumps with randomized magnitudes, and our Lyapunov function-based approach (see Appendices 1--3) is still valid.
The Poisson process captures the randomness of the arrival of platoons.

We are now ready to introduce the stochastic fluid model. 
\begin{dfn} \textbf{The stochastic fluid model} is defined as the tuple $\langle\mathcal Q,\mathcal U,G,\lambda,\mathcal V,S\rangle$, where
\begin{itemize}
    \item[-] $\mathcal Q:=\mathbb R_{\ge0}^3\times[0,\Theta]$ is the state space as well as the set of initial conditions,
    \item[-] $\mathcal U\subseteq\mathbb R_{\ge0}$ is the set of control inputs to the vector field $G$ (called the ``gate discharge''),
    \item[-] $G:\mathcal Q\times\mathcal U\to\mathbb R^4$ is a vector field such that $(d/dt)Q(t)=G(Q(t))$,
    \item[-] $\lambda\in\mathbb R_{\ge0}$ is the Poisson rate at which resets occur (called the ``arrival rate''),
    \item[-] $\mathcal V\subseteq\mathbb R^3$ is the set of control inputs to the reset mapping $S$ (called the ``allocation''),
    \item[-] $S:\mathcal Q\times\mathcal V\to\mathcal Q$ is a mapping that resets the state when a platoon arrives.
\end{itemize}
\end{dfn}

In our model, queuing happens due to (i) sudden increases in queues, which occur at rate $\lambda$ and according to the reset mapping $S$ , and (ii) interaction between queues in various links, which is captured by the vector field $G$.
We consider the control inputs $(u,v)$ to be determined by a \emph{control policy} $(\mu,\nu)$, i.e. $(u,v)=(\mu(q),\nu(q))$:

\begin{dfn}[Control policy]
    A control policy $(\mu,\nu)$ is specified by functions $\mu:\mathcal Q\to\mathcal U$ and $\nu:\mathcal Q\to\mathcal V$.
\end{dfn}
\noindent
We will describe how control policies for the stochastic fluid model can be translated to platoon coordination strategies in Section~\ref{sec_control}.

Given a control policy $(\mu,\nu)$, the model's dynamics can be expressed via the \emph{infinitesimal generator}
\begin{align}
    &\mathscr Lg(q)=G^T\Big(q;\mu(q)\Big)\nabla_qg(q)+\lambda\bigg(g\Big(S(q;\nu(q))\Big)-g(q)\bigg),\nonumber\\
    &\qquad q\in\mathcal Q,
\end{align}
where $g$ is a differentiable function \cite{davis84}. In the above, the first term on the right-hand side results from the fluid dynamics governed by the vector field $G$, and the second term results from the resets governed by the reset mapping $S$.

In summary, our model (as well as the subsequent analysis) focuses on the impact due to the following parameters:
\begin{enumerate}[(i)]
    \item Total demand $a$,
    \item Platooning ratio $\eta$ (or equivalently, platoon arrival rate $\lambda$),
    \item Platoon size $l$,
    \item Buffer size $\Theta$.
\end{enumerate}

The rest of this subsection is devoted to specifying the elements in the tuple $\langle\mathcal Q,\mathcal U,G,\lambda,\mathcal V,S\rangle$.

\subsubsection{State space $\mathcal Q$}

We use $q_0^m$ to denote the CAVs held in the gate, $q_1^m$ and $q_2^m$ to denote the queues of mainline traffic in links 1 and 2, respectively, and $q_1^o$ to denote the queue of off-ramp traffic in link 1.
The \emph{state} of the stochastic fluid model is $q=[q_0^m\ q_1^m\ q_1^o\ q_2^m]^T\in\mathcal Q$.
Note that $q_0^m$ consists of only CAVs, $q_1^m$ and $q_2^m$ consist of both CAVs and non-CAVs, and $q_1^o$ consists of only non-CAVs.
A key characteristic of platooning is the reduced inter-vehicle spacing.
\begin{figure}[hbt]
\centering
\includegraphics[width=0.2\textwidth]{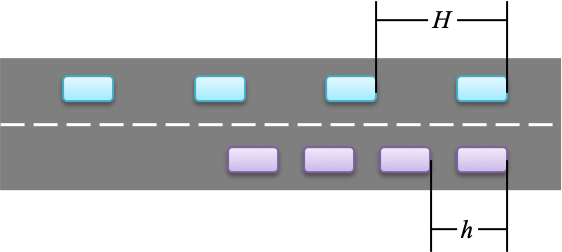}
\caption{$H$ is the spacing between ordinary vehicles, while $h$ is the spacing within a platoon; $\gamma=H/h$.}
\label{gamma}
\end{figure}
We model this by scaling down the CAV part of $q_k^m$ by a factor $\gamma>1$: a platoon of $l$ CAVs are roughly equivalent to $l/\gamma$ non-CAVs in terms of the occupied road space, where $\gamma$ is the ratio between the inter-non-CAV spacing $H$ and the inter-CAV spacing $h$; see Fig.~\ref{gamma}. A typical value for $\gamma$ is 2 \cite{bess+16procieee,santini2016consensus}.
Hence, $q$ is the vector of \emph{effective} queue lengths with the CAV part scaled down, which is in general smaller than the \emph{nominal} queue lengths.
Throughout this paper, we use $Q(t)=[Q_0^m(t)\ Q_1^m(t)\ Q_1^o(t)\ Q_2^m(t)]^T$ to denote the vector of queues at time $t$ and $q$ to denote a particular state.

\subsubsection{Gate discharge $\mathcal U$}

The rate at which the gate discharges traffic to link 1 is $u\in\mathcal U$.
We consider $u$ to be controlled by a \emph{gate discharge policy} $\mu:\mathcal Q\to\mathbb R_{\ge0}$ satisfying the following:
\begin{asm}
\label{asm_mu}
The gate discharge policy $\mu$ satisfies the following:
\begin{enumerate}[(i)]
    \item $\mu(q)$ is non-negative, bounded, and piecewise-continuous in $q$;
    \item $\mu(q)=0$ for $q$ such that $q_0^m=0$;
    \item $\mu(q)$ is non-increasing in $q_1^m,q_1^o,q_2^m$.
\end{enumerate}
\end{asm}
\noindent
In the above, (i) ensures regularity to facilitate analysis.
(ii) means that the gate discharge must vanish if $q_0^m=0$.
(iii) means that CAVs in the gate will be discharged slower if the downstream queues are longer.
We use $\mathscr U$ to denote the set of {gate discharge policies} satisfying the above assumption.
This assumption basically ensures that $Q(t)$ is bounded and piecewise-continuous in $t$.

\subsubsection{Vector field $G$}
\label{subsub_G}

\begin{figure}[hbt]
\centering
\includegraphics[width=0.35\textwidth]{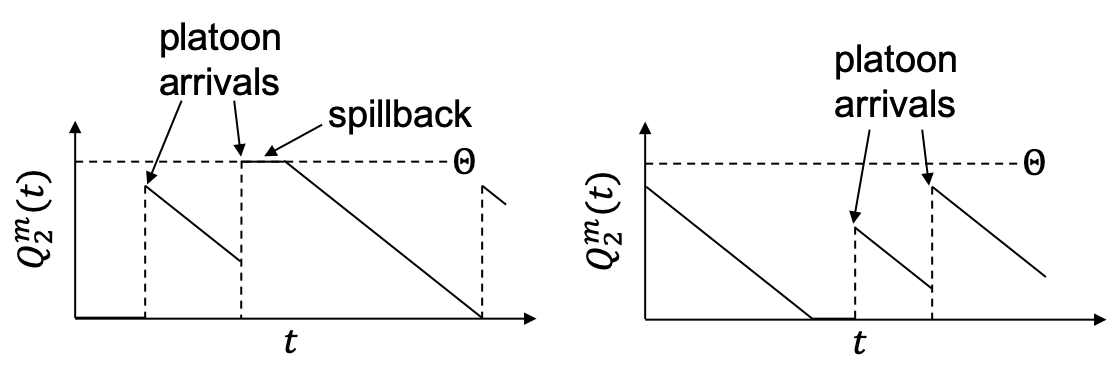}
\caption{Two examples of trajectories of $Q_2^m(t)$ with (left) and without spillback (right).}
\label{fig_spillback}
\end{figure}

The vector field $G$ specifies the model's dynamics between resets.
Specifically, the inflow of mainline background traffic is $(1-\eta)\rho a$, and the inflow of off-ramp background traffic is $(1-\rho)a$.
If the queue in link 2 is less than the buffer size, then the queue in link 1 is discharged at the link's capacity $F$; otherwise, the queue propagates to link 1 and reduces the off-ramp flow (called \emph{spillback}, see Fig.~\ref{fig_spillback}). 
The flows in Fig.~\ref{fig_tandem} are given by
\begin{subequations}
\begin{align}
&f_0(q;\mu(q)):=\mu(q),\\
&f_1(q;\mu(q)):=\nonumber\\
&\left\{\begin{array}{ll}
\min\{(1-\eta)\rho a+\mu(q),F\} & \mbox{if $q_1^m=0,q_2^m<\Theta$}, \\
F & \mbox{if $q_1^m>0,q_2^m<\Theta$}, \\
\min\{(1-\eta)\rho a+\mu(q),F-R\} & \mbox{if $q_1^m=0,q_2^m=\Theta$}, \\
F-R & \mbox{if $q_1^m>0,q_2^m=\Theta$},
\end{array}\right.\label{eq_f1}\\
&f_{2}(q;\mu(q)):=\left\{\begin{array}{ll}
\min\{f_1(q;\mu(q)),F-R\} & \mbox{if $q_2^m=0$}, \\
F-R & \mbox{if $q_2^m>0$}
\end{array}\right.\\
&{r}(q;\mu(q)):=\nonumber\\
&{\footnotesize\left\{\begin{array}{ll}
\min\{(1-\rho)a,F-f_1(q;\mu(q)),R\} & \mbox{if $q_1^o=0,q_2^m<\Theta$}, \\
\min\{F-f_1(q;\mu(q)),R\} & \mbox{if $q_1^o>0,q_2^m<\Theta$}, \\
\min\{(1-\rho)a,F-R-f_1(q;\mu(q)),R\} & \mbox{if $q_1^o=0,q_2^m=\Theta$}, \\
\min\{F-R-f_1(q;\mu(q)),R\} & \mbox{if $q_1^o>0,q_2^m=\Theta$}.
\end{array}\right.}\label{eq_r}
\end{align}
\end{subequations}
\eqref{eq_f1} and \eqref{eq_r} indicate that spillback happens if $q_2^m=\Theta$: whenever the threshold is attained, the upstream capacity is dropped.
To focus on the impact of capacity, we assume that traffic queues are always discharged at the rate of $F$ or $F-R$.\footnote{Simulation results show that this simplified flow model is largely consistent with more sophisticated models such as the CTM; see Section Section~\ref{sub_ctm}.}
Then, the fluid dynamics is specified by the vector field $G:\mathcal Q\times\mathcal U\to\mathbb R^4$ defined as
\begin{subequations}
\begin{align}
&G_0^m(q;\mu(q)):=-\mu(q),\label{eq_G0}\\
&G_1^m(q;\mu(q)):=(1-\eta)\rho a-(f_1(q;\mu(q))+{r}(q;\mu(q))),\\
&G_1^o(q;\mu(q)):=(1-\rho)a-{r}(q;\mu(q)),\\
&G_2^m(q;\mu(q)):=f_1(q;\mu(q))-f_2(q;\mu(q)).\label{eq_G2}
\end{align}
\end{subequations}
To emphasize that the vector field depends on the control policy, we use the notation $G(q,\mu(q))$.
Note that for an admissible gate discharge policy $\mu\in\mathscr U$, $G$ is bounded and piecewise-continuous in $q$.
Since we focus on the aggregate behavior of both traffic classes, our results hold for a variety of capacity-sharing models (see e.g. \cite{jin2019behavior}). In addition, microscopic maneuvers such as overtaking are implicitly captured by the flow dynamics.

\subsubsection{Arrival rate $\lambda$}
CAV platoons arrive according to a Poisson process of rate $\lambda$ (per hr), which is given by \eqref{eq_lambda}.
$\lambda$ specifies the rate at which the continuous state $Q(t)$ is reset.
Randomness of the platoon arrival process can be attributed to the process of platoon formation \cite{xiong2019analysis,xiong2020optimizing}.

\subsubsection{Allocation $\mathcal V$}
Platoon control is modeled by a vector $v=[v_0\ v_1\ v_2]^T\in\mathcal V$. Arriving platoons are allocated to each link according to $v$: for example, $v=[l,0,0]^T$ means that a platoon is allocated to the gate. Recall that $v$ is determined by a mapping $\nu:\mathcal Q\to\mathcal V$.
We assume the following for $\nu$:
\begin{asm}
\label{asm_nu}
The allocation policy $\nu$ satisfies the following:
\begin{enumerate}[(i)]
    \item $\sum_{k=0}^2\nu_k(q)=0$, $0\le\nu_0(q)\le l/\gamma$, $-l/\gamma\le\nu_1(q)\le0$, $\max\{-l/\gamma,q_2^m-\Theta\}\le\nu_2(q)\le0$,
    \item $\nu_2(q)=0$ for $q$ such that $q_1^m>0$, and $\nu_1(q)=\nu_2(q)=0$ for $q$ such that $q_0^m>0$,
    \item $\nu_k(q)$ is non-increasing in $q_k^m$ and non-decreasing in $q_j^m$ for $j\neq k$, and $\nu_k(q)$ is non-increasing in $q_1^o$ for $k=1$ and non-decreasing in $q_1^o$ for $k\neq 1$.
\end{enumerate}
\end{asm}
\noindent
In the above, (i) means that $\nu$ only {distributes} but does not {creates} traffic. 
(ii) results from the ``first-come-first-served'' principle: a platoon cannot be allocated to link $k$ if there is a non-zero queue in link $k-1$.
(iii) means that more traffic is allocated to a link with a shorter queue.
We use $\mathscr V$ to denote the set of {gate discharge policies} satisfying the above assumption.

\subsubsection{Reset mapping $S$}
\label{subsub_S}

Arrivals of platoons lead to sudden increases in the state $Q(t)$, and the reset mapping is given by
\begin{subequations}
\begin{align}
&S_0^m(q;v):= q_0^m+v_0,\label{eq_T0}\\
&S_1^m(q;v):= q_1^m+(q_2^m+l/\gamma-\Theta)_++v_1,\label{eq_T1}\\
&S_1^o(q;v):= q_1^o,\\
&S_2^m(q;v):= \min\{\Theta,q_2^m+l/\gamma\}+v_2.\label{eq_T2}
\end{align}
\end{subequations}
In particular, $S(q;0)$ represents the reset mapping if no control is applied (i.e. platoons are not coordinated).\footnote{If $v=0$, then no CAVs will ever be allocated to the gate; hence the gate discharge has no impact.}
For $v=0$, no platoons will be allocated to the gate upon arrival; instead, every platoon will be allocated to link 2 unless $Q_2^m(t)$ attains the buffer size $\Theta$.
If a platoon arrives at time $t$, then the state is reset according to
\begin{align*}
    Q(t)=S\Big(Q(t_-);v\Big),
\end{align*}
where $Q(t_-)$ is the vector of queues immediately before the arrival.


\subsection{Problem definition}
\label{sub_preliminaries}

The main questions that we study are
\begin{enumerate}[(i)]
    \item Given model parameters and a control policy, how to determine whether the queues are bounded (in expectation), and how to compute or estimate the model's throughput?
    \item How to design the control policy to ensure bounded queues and to improve throughput and travel time?
\end{enumerate}
To study the above questions, we introduce the following definitions.

First, following \cite{dai1995stability}, we define stability as follows:
\begin{dfn}[Stability]
The stochastic fluid model is {stable} if there exists $Z<\infty$ such that for each initial condition $q\in\mathcal Q$
\begin{align}\label{eq_bounded}
\limsup_{t\to\infty}\frac1t\int_{s=0}^t{\mathrm E}\left[|Q(s)|\right]ds\le Z.
\end{align}
\end{dfn}
\noindent 
Practically, stability means that expected queue size is bounded, and hence the probability of long queues is small.

Second, given a control policy $(\mu,\nu)$, the stochastic fluid model typically admits an \emph{invariant set}, which is defined as follows:
\begin{dfn}[Invariant set]
Given a control policy $(\mu,\nu)$, a compact set $\mathcal M_{\mu,\nu}\subseteq\mathcal Q$ is an invariant set if
\begin{enumerate}[(i)]
    \item $\lim_{t\to\infty}\Pr\{Q(t)\in\mathcal M_{\mu,\nu}|Q(0)=q\}=1,
    \quad\forall q\in\mathcal Q,$
    \item $Q(t)\in\mathcal M_{\mu,\nu}, \quad\forall Q(0)=q\in\mathcal M_{\mu,\nu}.$
\end{enumerate}
\end{dfn}
\noindent The interpretation of an invariant set is that (i) for each initial condition, the process $\{Q(t);t\ge0\}$ enters the set $\mathcal M_{\mu,\nu}$ almost surely (a.s.), and (ii) if the process $\{Q(t);t\ge0\}$ starts within $\mathcal M_{\mu,\nu}$, then it never leaves $\mathcal M_{\mu,\nu}$.
Since stability is defined for an infinite time horizon in \eqref{eq_bounded}, we can focus on the model's evolution over $\mathcal M_{\mu,\nu}$ rather than over $\mathcal Q$; this simplifies the analysis.
Also note that $\mathcal M_{\mu,\nu}$ depends on $(\mu,\nu)$ and is thus typically determined based on characteristics of $(\mu,\nu)$.

The theoretical tool that we use to establish stability is the \emph{Foster-Lyapunov criterion}, which is a sufficient condition for \eqref{eq_bounded}:

\noindent{\bf Foster-Lyapunov criterion \cite{meyn93}}.
\emph{Consider a Markov process with an invariant set $\mathcal Y$ and infinitesimal generator $\mathscr A$. If there exist a Lyapunov function $W:\mathcal Y\to\mathbb R_{\ge0}$ and constants $c>0$, $d<\infty$ satisfying
\begin{align}
    \mathscr AW(y)\le-cg(y)+d,
    \quad\forall y\in\mathcal Y,
    \label{eq_AW}
\end{align}
then for each initial condition $y\in\mathcal Y$
\begin{align}
    \limsup_{t\to\infty}\frac1t\int_{\tau=0}^t\mathrm E[g(Y(t))]d\tau\le d/c.
\end{align}
}%
In the above, \eqref{eq_AW} is called the ``drift condition'' \cite{meyn93}.
Verifying the drift condition is in general challenging, since it requires (i) finding an effective Lyapunov function, which is not straightforward for a nonlinear system as our stochastic fluid model, and (ii) checking the inequality \eqref{eq_AW} over a possibly unbounded set $\mathcal Y$, which involves a non-convex optimization.
In Section~\ref{sec_analysis}, we argue how we address these challenges.

Third, with the notion of stability, we define the stochastic fluid model's \emph{throughput} as follows:
\begin{dfn}[Throughput]
Given a control policy $(\mu,\nu)\in\mathscr U\times\mathscr V$, the throughput of the stochastic fluid model is
\begin{align*}
    \bar a_{\mu,\nu}:=&\sup a\\
&\mbox{s.t.} \mbox{ fluid model is stable under }(\mu,\nu)\in\mathscr U\times\mathscr V.
\end{align*}
\end{dfn}
\noindent 
$\bar a_{\mu,\nu}$ is defined as supremum rather than maximum, since the stability constraint may lead to strict inequalities.
As indicated in the above definition, the key to throughput analysis is to develop stability conditions for the stochastic fluid model, which we discuss in Section~\ref{sec_analysis}.
We have the following preliminary result for throughput:
\begin{lmm}[Nominal throughput]
\label{lmm_maximum}
For any control policy $(\mu,\nu)\in\mathscr U\times\mathscr V$, throughput $\bar a_{\mu,\nu}$ of the stochastic fluid model is upper-bounded by
\begin{align}
\bar a_{\mu,\nu}\le a^*:=\min\Big\{
\frac{R}{1-\rho},
\frac{F-R}{(\eta/\gamma+1-\eta)\rho}
\Big\}.
\label{eq_nominal}
\end{align}
\end{lmm}
We call $a^*$ as defined in \eqref{eq_nominal} the \emph{nominal throughput}.
To interpret the expression for $a^*$, note that the first (resp. second) term in $\min\{\cdot\}$ results from the capacity constraint of the off-ramp (resp. the mainline bottleneck).
Importantly, we will show that the nominal throughput cannot always be attained due to the interaction between CAV and non-CAV traffic and due to spillback of traffic queues. 

Finally, for control design, we consider the following formulation:
\begin{align*}
\mbox{decision:}&\quad (\mu,\nu)\in\mathscr U\times\mathscr V\nonumber\\
\mbox{(P)}\quad\mbox{objective:}&\quad\max\quad \bar a_{\mu,\nu}\\
&\quad s.t.\quad \mbox{stability}.
\end{align*}
That is, the first objective is stabilization, and the second objective is throughput maximization or queue minimization.
Note that the objective value of (P) is upper-bounded by $a^*$ in Lemma~\ref{lmm_maximum}.

%% file: Texts/3_Analysis.tex
\section{Stability and throughput}
\label{sec_analysis}

In this section, we study the stability and throughput of the stochastic fluid model, which provides insights about the efficiency of the mixed-traffic highway.

The first main result gives a criterion to check the stability of the stochastic fluid model under a given control policy $(\mu,\nu)$.
\begin{thm}[Stability criterion]
\label{thm_sufficient}
Suppose that the stochastic fluid model admits an invariant set $\mathcal M_{\mu,\nu}\subseteq\mathcal Q$.
The stochastic fluid model is stable if
\begin{subequations}
\begin{align}
    & a<a^*,\quad\mbox{and}
    \label{eq_a<a*}\\
    &\max_{\substack{\xi\in\mathcal M_{\mu,\nu}:\\\xi_0^m=\xi_1^m=\xi_1^o=0}}
    \frac{\xi_2^m}{\Theta}\Big((1-\eta)\rho a-(F-R)\Big)\nonumber\\
    &\hspace{1.5cm}+\frac{\eta\rho a}{l} \bigg(S^m_0(\xi;\nu(\xi))+S^m_1(\xi;\nu(\xi))\nonumber\\
    &\hspace{1.5cm}+\frac{1}{2\Theta}\Big(S^m_2(\xi;\nu(\xi))\Big)^2-\frac1{2\Theta}(\xi_2^m)^2\bigg)\nonumber\\
    &\quad<\frac{R-(1-\rho)a}{(1-\rho)a}\Big((F-R)-(\eta/\gamma+1-\eta)\rho a\Big).\label{eq_max}
\end{align}
\end{subequations}
\end{thm}

One can interpret the stability criterion as follows.
\eqref{eq_a<a*} results from the nominal upper bound in Lemma~\ref{lmm_maximum}.
\eqref{eq_max} essentially results from the interaction between the mainline and the off-ramp traffic.
\eqref{eq_max} also captures the influence of the control policy $(\mu,\nu)$.
Although the complexity of the maximization on the left-hand side of \eqref{eq_max} depends on the control policy $(\mu,\nu)$, the decision variable $\xi$ can only vary in the direction of $\xi_2^m$; all the other components of $\xi$ must be zero.
Hence, the maximization involves essentially only one decision variable $\xi_2^m$, which takes values from a compact interval $[0,\Theta]$, and is thus not hard to solve numerically.
This is a significant refinement of the Foster-Lyapunov criterion, which, in its general form, does not give a ready-to-use stability criterion for our stochastic fluid model.
Theorem~\ref{thm_sufficient} can be used for throughput analysis by finding the largest demand $a$ that satisfies \eqref{eq_a<a*}--\eqref{eq_max}.
Since Theorem~\ref{thm_sufficient} is a sufficient condition for stability, it leads to a lower bound for throughput.

We prove Theorem~\ref{thm_sufficient} by considering a quadratic Lyapunov function and establishing the Foster-Lyapunov criterion \cite{meyn1993stability} for the queuing process. The main technique is to relate the fluid queuing process to an underlying M/D/1 process, which we discuss below. The detailed proof is in Appendix 1.
\vspace{-12pt}
\begin{figure}[hbt]
\centering
\subfigure[Trajectories.]{
\centering
\includegraphics[width=0.3\textwidth]{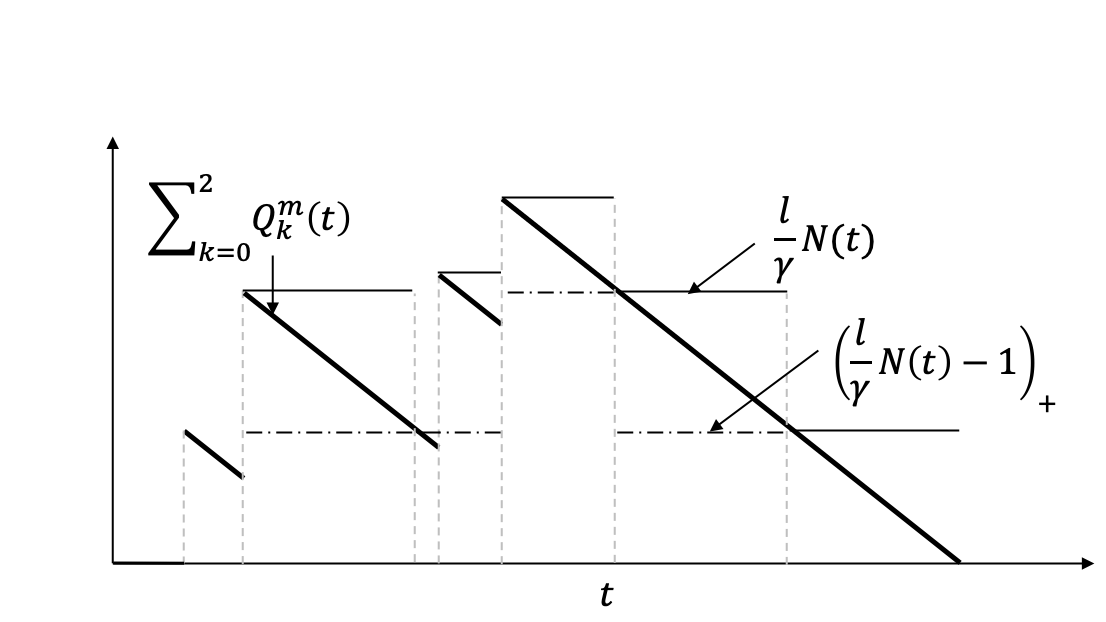}
\label{md1}
}
\subfigure[Steady-state CDFs.]{
\centering
\includegraphics[trim=0 2cm 0 0,clip,width=0.15\textwidth]{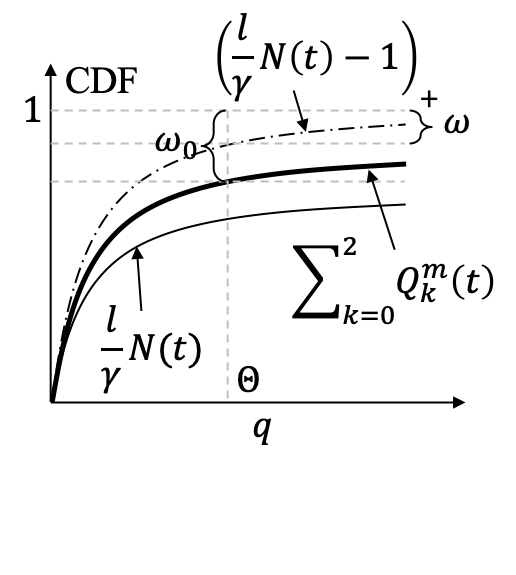}
\label{cdf}
}
\caption{Relation between M/D/1 process $N(t)$ and fluid process $\sum_{k=0}^2Q_k^m(t)$.}
\end{figure}
%

\begin{figure*}[hbt]
\centering
\subfigure[Throughput vs. fraction.]{
\centering
\includegraphics[width=0.31\textwidth]{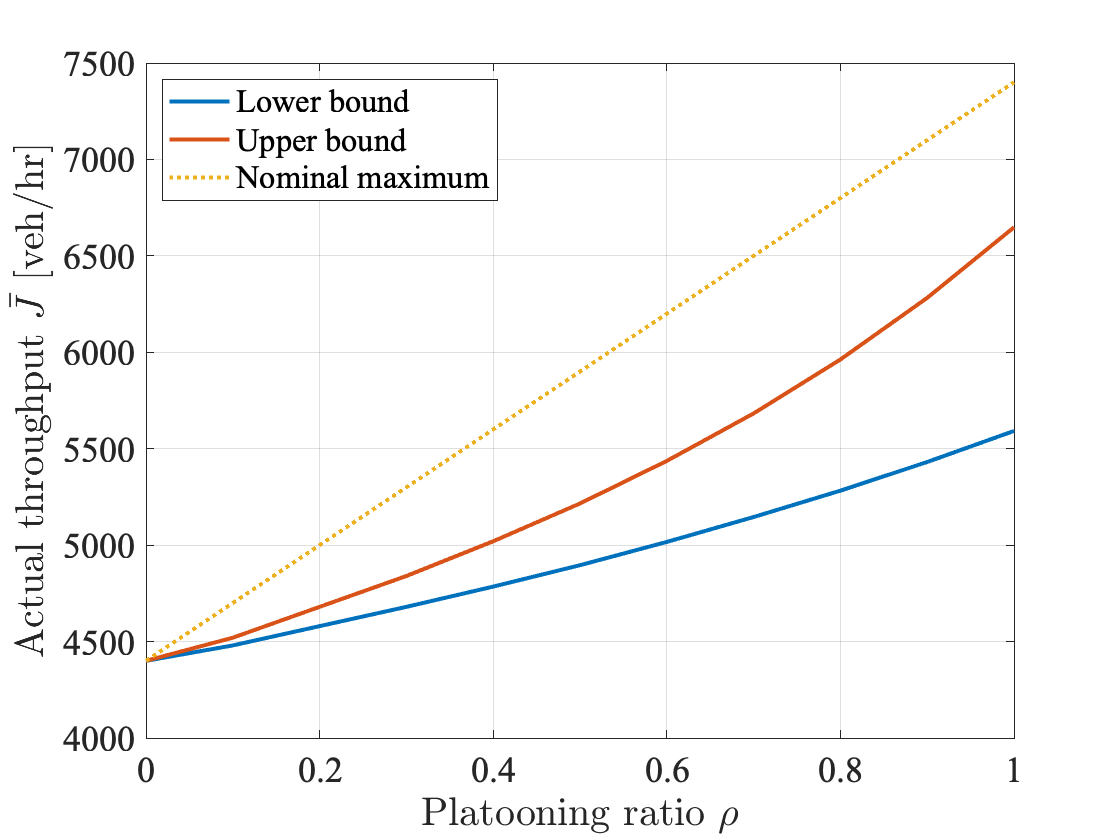}
\label{Rbar_eta}
}
\subfigure[Throughput vs. platoon size.]{
\centering
\includegraphics[width=.31\textwidth]{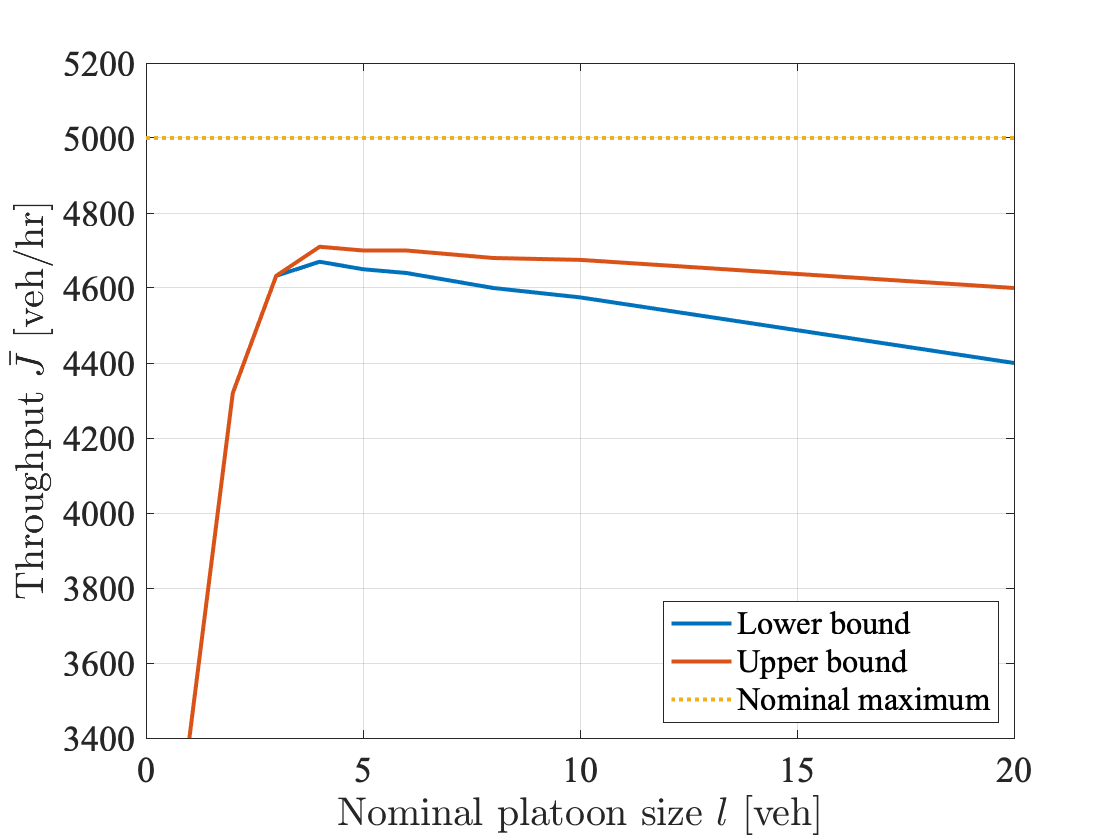}
\label{Rbar_l}
}
\subfigure[Throughput vs. buffer size.]{
\centering
\includegraphics[width=0.31\textwidth]{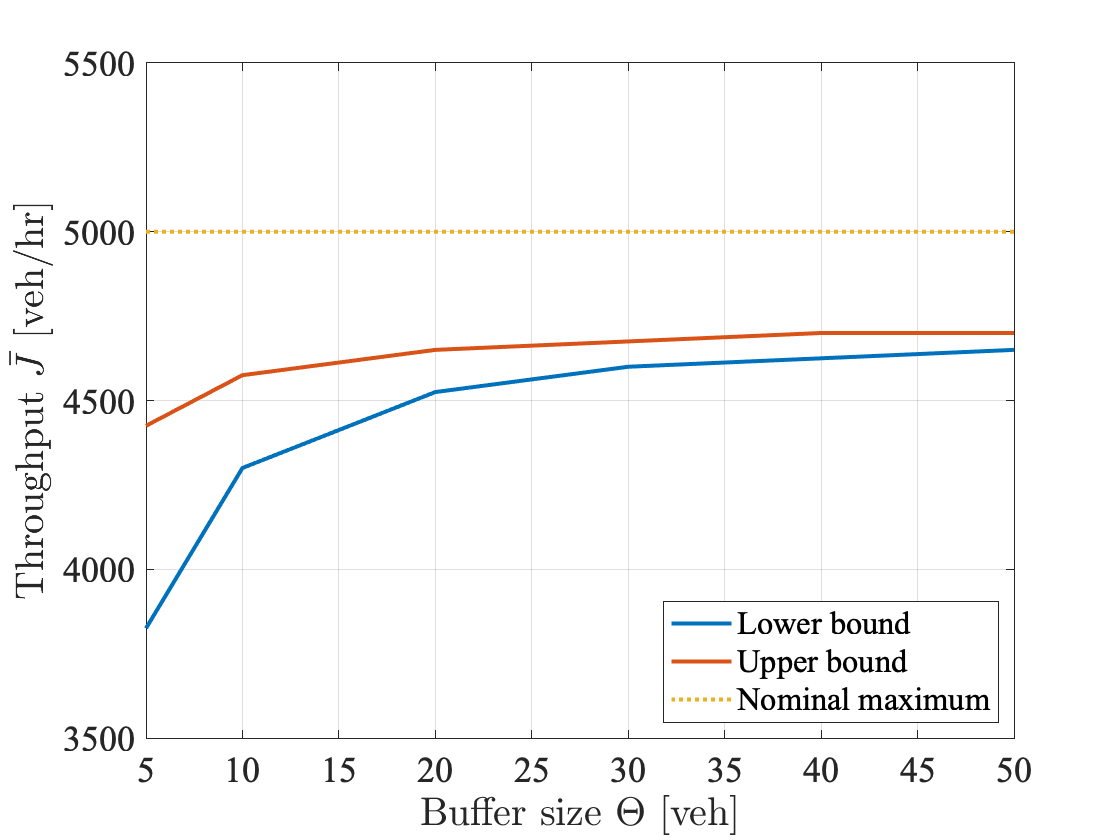}
\label{fig_bounds}
}
\caption{Relation between throughput and model parameters. Bounds are computed using Theorem~\ref{thm_bounds}. ``Nominal maximum'' refers to $a^*$ in Lemma~\ref{lmm_maximum}. In Fig.~\ref{Rbar_eta} (resp. Fig.~\ref{Rbar_l}), $\eta=0$ (resp. $l=1$) means no platooning.}
\end{figure*}

In particular, if no control is applied, then \eqref{eq_max} can be manually solved, and explicit lower and upper bounds for throughput can be derived. 
To see this, consider the process
\begin{align}\label{eq_N}
N(t):=\Big\lceil\frac{\gamma}l\Big(\sum_{k=0}^2Q_k^m(t)\Big)\Big\rceil
\quad t\ge0
\end{align}
which satisfies the following:
\begin{lmm}
The process $\{N(t);t\ge0\}$ is an M/D/1 process with arrival rate $\lambda=\eta\rho a/l$ and service time $s:=\frac{l}{\gamma(F-R-(1-\eta)\rho a))}$. Furthermore, 
\begin{align}
    \Big((l/\gamma)N(t)-1\Big)_+\le \sum_{k=0}^2Q_k^m(t)\le(l/\gamma)N(t),
    \quad \forall t\ge0.
\end{align}
\end{lmm}
\noindent
That is, we can use $N(t)$ to bound $\sum_{k=0}^2Q_k^m(t)$; see Fig.~\ref{md1}.
The steady-state probabilities $\pi_n$ of the M/D/1 process are given by a standard result in queuing theory \cite{shortle2018fundamentals}:
\begin{align}
&\pi_n=\left\{\begin{array}{l}
1-\lambda s 
\quad n=0,\\
(1-\lambda s)(e^{\lambda s} - 1) 
\quad n=1,\\
(1-\lambda s)\Big(e^{n\lambda s}
+\sum_{k=1}^{n-1}e^{k\lambda s}(-1)^{n-k}\\
\quad\times\left[\frac{(k\lambda s)^{n-k}}{(n-k)!}+\frac{(k\lambda s)^{n-k-1}}{(n-k-1)!}\right]\Big)
\quad n\ge2.
\end{array}\right.\label{eq_pin}
\end{align}
Using the above result, we can obtain a lower bound
\begin{align}
    \omega=1-\sum_{n=0}^{\lfloor{\gamma\Theta/l}\rceil}\pi_n,\label{eq_omega}
\end{align}
for the actual fraction of time $\omega_0$ that the stochastic fluid model experiences spillback; see Fig.~\ref{cdf}.

With the above arguments, we can state the second main result of this section as follows; the proof is in Appendix 2.

\begin{thm}[Throughput without control]
\label{thm_bounds}
Suppose that the stochastic fluid model is not controlled, i.e. $\mu(q)=0$ and $\nu(q)=0$ for all $q\in\mathcal Q$.
Then, the throughput $\bar a$ of the model is bounded by
\begin{align}
&\min\Big\{
\frac{F-R}{\rho(\eta/\gamma+1-\eta)},
\frac{R}{1-\rho+\frac{1}{2}\Big(\sqrt{\zeta^2+\frac{2\rho Rl}{\gamma\Theta(F-R)}}-\zeta)\Big)}
\Big\}\nonumber\\
&\quad
\le
\bar a\le
\min\Big\{
\frac{F-R}{\rho(\eta/\gamma+1-\eta)},
\frac{(1-\omega)R}{(1-\rho)}\Big\}
\end{align}
where
\begin{align}
\zeta=(1-\rho)-\rho(\eta/\gamma+1-\eta)\frac{R}{F-R}\label{eq_zeta0}
\end{align}
and $\omega$ is given by \eqref{eq_omega}.
\end{thm}

Using the bounds in Theorem~\ref{thm_bounds}, we can analyze the model's throughput without control.
Figs.~\ref{Rbar_eta}--\ref{fig_bounds} summarize our results for throughput analysis with the nominal parameters in Table~\ref{tab_parameters}.
Importantly, due to the interaction between multiple traffic classes and due to lack of coordination between platoons, the nominal throughput given by Lemma~\ref{lmm_maximum} is not attained.
Specific discussions about the figures are as follows:



%
\begin{table}[h]
\caption{Nominal model parameter values.}
\label{tab_parameters}
\centering
\begin{tabular}{@{}lcc@{}}
\toprule
Quantity                   & Notation & Nominal value       \\ \midrule
mainline capacity          & $F$       & 4500 veh/hr \\ \hline
ramp capacity          & $R$       & 1500 veh/hr \\ \hline
platooning ratio & $\eta$        & 0.2 \\ \hline
platoon scaling coefficient & $\gamma$        & 2 \\ \hline
platoon size             & $l$        & 5 veh       \\ \hline
mainline ratio & $\rho$        & 0.75 \\ \hline
buffer size & $\Theta$        & 50 veh \\ \bottomrule
\end{tabular}
\end{table}

\begin{enumerate}
    \item Fig.~\ref{Rbar_eta} shows that higher fraction of platooning lead to higher throughput, which is consistent with previous results \cite{lioris17}. The nominal throughput is not attained due to interaction between CAVs and non-CAVs.

    \item Fig.~\ref{Rbar_l} shows that platooning does improve throughput, but overly long platoons may result in lower throughput.
The reason is that longer platoons have stronger impact on local traffic and cause larger local congestion. For this example, the empirically optimal platoon size is 4 CAVs.

        \item As shown in Fig.~\ref{fig_bounds}, when the buffer size $\Theta$ is small (e.g. less than 20), spillback occurs frequently, and thus an obvious throughput drop (with respect to the nominal throughput) is observed.
As $\Theta$ approaches infinity, spillback hardly occurs, and both bounds approach the nominal throughput given by Lemma~\ref{lmm_maximum}.
\end{enumerate}

%% file: Texts/4_Control.tex
\section{Control Design}
\label{sec_control}

In this section, we study a set of control policies that attain the nominal throughput given by Lemma~\ref{lmm_maximum}.

One can indeed use Theorem~\eqref{thm_sufficient} as the stability constraint and solve (P).
However, since Theorem~\ref{thm_sufficient} is a sufficient condition, it  in general leads to a sub-optimal solution.
The main result of this section is a sufficient condition for optimality of a given control policy:

\begin{thm}[Optimality criterion]
\label{thm_optimal}
Suppose that a control policy $(\mu,\nu)\in\mathscr U\times\mathscr V$ satisfies
\begin{subequations}
\begin{align}
    &\mu(q)=0,
    \quad\forall q:q_1^m>0\mbox{ or } q_2^m=\Theta,\label{eq_muq=0}\\
    &\mu(q)\le F-((1-\eta)\rho+(1-\rho))a,
    \quad\forall q\in\mathcal Q,\label{eq_muq<F}\\
    &\mu(q)+(1-\eta)\rho a\ge F-R,
    \quad\forall q:q_0^m>0\mbox{ and }q_2^m=0,\label{eq_muq>F-R}\\
    &S_1^m(q;\nu(q))=0,\quad\forall q\in\mathcal Q,\label{eq_S1m}\\
    &S_2^m(q;\nu(q))<\Theta,
    \quad\forall q\in\mathcal Q.
    \label{eq_S2m}
\end{align}
\end{subequations}
Then, $(\mu,\nu)$ stabilizes the stochastic fluid model if and only if
\begin{align}
    a<a^*,\label{eq_iff}
\end{align}
where $a^*$ is the nominal throughput given by \eqref{eq_nominal}.
Furthermore, under \eqref{eq_iff}, the control $(\mu,\nu)$ minimizes the total queue length $|Q(t)|$ among all admissible controls at any time $t\ge0$, and the time-average queuing delay converges as follows:
\begin{align}
&\lim_{t\to\infty}\frac1t\int_{\tau=0}^t|Q(\tau)|d\tau
\stackrel{\footnotesize a.s.}=\bar Q
:=\frac{\eta\rho al}{2\gamma^2(F-R-(1-\eta)\rho a)}\nonumber\\
&\qquad\times\bigg(\frac{\eta\rho a}{\gamma(F-R-(\eta/\gamma+1-\eta)\rho a)}+1\bigg).
\label{eq_Qbar}
\end{align}
\end{thm}

Note that Theorem 3 addresses not only stabilization, but also throughput maximization and queue minimization.
The optimality criterion \eqref{eq_muq=0}--\eqref{eq_S1m} can be interpreted as follows. \eqref{eq_muq=0}, \eqref{eq_muq<F}, and \eqref{eq_S1m} ensure no queuing in link 1 (i.e. $Q_1^m(t)=0$).
\eqref{eq_muq=0} and \eqref{eq_S2m} ensure no spillback at link 2 (i.e. $Q_2^m(t)<\Theta$).
\eqref{eq_muq>F-R} ensures that as long as there is a non-zero queue in the gate (i.e. $Q_0^m(t)>0$), link 2 must be discharging traffic at its capacity $F-R$.

We use $\mathscr U^*\times\mathscr V^*$ to denote the set of control policies satisfying \eqref{eq_muq=0}--\eqref{eq_S2m}.
Since $(\mu,\nu)\in\mathscr U^*\times\mathscr V^*$ stabilizes the stochastic fluid model if and only if the demand is less than the nominal throughput, $(\mu,\nu)$ maximizes the throughput.
In addition, each control in $\mathscr U^*\times\mathscr V^*$ not only maximizes throughput, but also minimizes queuing delay in a sample path-wise manner.
Furthermore, an analytical expression for the mean queuing delay is obtained.
The proof of Theorem~\ref{thm_optimal} is in Appendix 3.



In the rest of this section, we discuss two concrete control policies with practical interpretations, viz. headway regulation and platoon size management.
Here we focus on their formulation in the stochastic fluid model.
In Section~\ref{sub_implementation}, we discuss and demonstrate how they can be implemented in practice (e.g. speed).

\subsection{Headway regulation}
\label{exm_ustar}
Under this strategy, if a large number of platoons arrive within a short time period, some platoons will be allocated to the gate so that their arrival at the bottleneck is postponed to avoid cumulative congestion at the bottleneck.
In the stochastic fluid model, the headway regulation strategy can be formulated as a control policy $(\mu,\nu)$ such that:
\begin{subequations}
\begin{align}
&\mu^{hr}(q):=\left\{\begin{array}{ll}
\alpha & \mbox{if }q_0^m>0,q_2^{m}=0,\\
\alpha & \mbox{if }\gamma q_0^m/l\notin\{0,1,2,\ldots\},\\
0 & \mbox{o.w.}
\end{array}\right.
\label{eq_ustar}\\
&\nu^{hr}(q):=\left[\begin{array}{c}
     l/\gamma  \\
     q_1^m-S_1(q;0)\\
     q_2^m-S_2(q;0)
\end{array}
\right],
\end{align}
\end{subequations}
where $\alpha$ is the saturation flow rate of CAVs given by
\begin{align}
\alpha=v_0/h,
\label{eq_alpha}
\end{align}
$v_0$ is the nominal speed of the highway section, and $h$ is the between-vehicle spacing within platoons.
One can check that $(\mu^{hr},\nu^{hr})\in\mathscr U^*\times\mathscr V^*$.

\begin{figure}[hbt]
\centering
\subfigure[Regulating headway.]{
\centering
\includegraphics[width=0.2\textwidth]{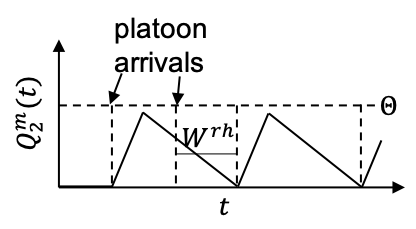}
\label{fig_Qt_controlled}
}
\subfigure[Splitting platoon.]{
\centering
\includegraphics[width=0.2\textwidth]{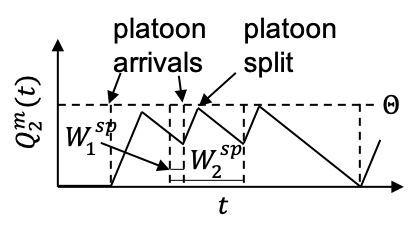}
\label{fig_Qt_controlled2}
}
\caption{Illustration of how regulating headway $(\mu^{hr},\nu^{hr})$ and splitting platoons $(\mu^{sm},\nu^{sm})$ avoid spillback; see Fig.~\ref{fig_spillback} for the uncontrolled case.}
\label{fig_Qt}
\end{figure}

Fig.~\ref{fig_Qt_controlled} illustrates the idea of this control policy: if two platoons enter the highway with a short inter-arrival time, the following platoon is decelerated so that its arrival at the bottleneck is postponed by $W^{hr}$ amount of time.
Consequently, the platoons arrive at the bottleneck with sufficient headway in between.
As illustrated in Fig.~\ref{fig_Qt_controlled}, $(\mu^{hr},\nu^{hr})$ essentially regulates times at which platoons arrive at the bottleneck so that congestion does not build up or spill back from the bottleneck, and the off-ramp traffic is not blocked.
To compute $W^{hr}$, suppose that a platoon enter the highway at time $t$ and the state immediately before the arrival is $Q(t_-)$; then $W^{hr}$ is the solution to the deterministic equation
\begin{align}
    \int_{s=t}^{W^{hr}}\mu^{hr}(Q(s))ds=Q_0^m(t_-).
    \label{eq_Whr}
\end{align}
Note that $W^{hr}$ is independent of any platoon arrivals after $t$.

\subsection{Platoon size management}
\label{exm_u2star}
If the highway is congested, long platoons will be disadvantageous at bottlenecks due to their sizes. Consequently, the operator can instruct platoons to split into shorter platoons to mitigate local congestion. The decision variable is whether to split or maintain a platoon as it enters the highway.
In the stochastic fluid model, splitting a platoon can be modeled by a control policy $(\mu^{sm},\nu^{sm})$ defined as follows:
\begin{align}
&{{\mu}^{sm}}(q):=\left\{\begin{array}{ll}
\alpha & \mbox{if }q_0^m>0,q_2^{a}\le \Theta-\frac{l((\eta/\gamma+1-\eta)\rho a-(F-R))}{2\gamma \eta\rho a},\\
\alpha & \mbox{if }(\gamma q^m_0/l)\notin\{0,\frac12,1,\frac32,\ldots\},\\
0 & \mbox{o.w.}
\end{array}\right.
\label{eq_u2star}\\
&\nu^{sm}(q):=\left[\begin{array}{c}
     l/\gamma  \\
     q_1^m-S_1(q;0)\\
     q_2^m-S_2(q;0)
\end{array}
\right].
\end{align}
One can check that $(\mu^{sm},\nu^{sm})\in\mathscr U^*\times\mathscr V^*$.
As shown in Fig.~\ref{fig_Qt_controlled2}, ${\mu}^{sm}$ opens the gate only if link 2 has sufficient space to accept at least half a platoon.
Furthermore, if $q_2^{a}$ is close to the buffer size $\Theta$, ${\nu}^{sm}$ will split a platoon into two short platoons and allocate the two short platoons to links 0 and 2, respectively, to avoid spillback.
In practice, suppose a platoon enters the highway at time $t$ and let $T$ be the solution to 
\begin{align}
&\int_{s=t}^Tf_2(Q(s))ds=Q^m_0(t_-)+Q^m_1(t_-)\nonumber\\
&\qquad\qquad+\min\Big\{\frac{l}{2\gamma},\Big(\Theta-\frac{l}{2\gamma}- Q_2^m(t_-)\Big)_+\Big\},
\label{eq_T}
\end{align}
where $Q_k^m(t_-)$ is the queue size in the gate immediately before the platoon arrives.
Then, the delays $W^{sm}_1$ and $W^{sm}_2$ indicated in Fig.~\ref{fig_Qt_controlled2} are given by
\begin{align}
W_1^{sm}=T,\
W_2^{sm}=T+\frac{l}{2\gamma(F-R-(1-\eta)\rho a)}.
\label{eq_W1}
\end{align}

%% file: Texts/4.5_Simulate.tex
\section{Implementation, simulation, and discussion}
\label{sec_ctm}

In this section, we translate the control laws discussed in the previous section to platoon coordination instructions that can be implemented in practice (Section~\ref{sub_implementation}).
We also validate the optimality of the headway regulation strategy via two standard simulation environments, viz. CTM and SUMO (Section~\ref{sub_simulation}).

\subsection{Implementation of proposed control policies}
\label{sub_implementation}

\begin{figure}[hbt]
\centering
\subfigure[Headway regulation.]{
\centering
\includegraphics[width=0.4\textwidth]{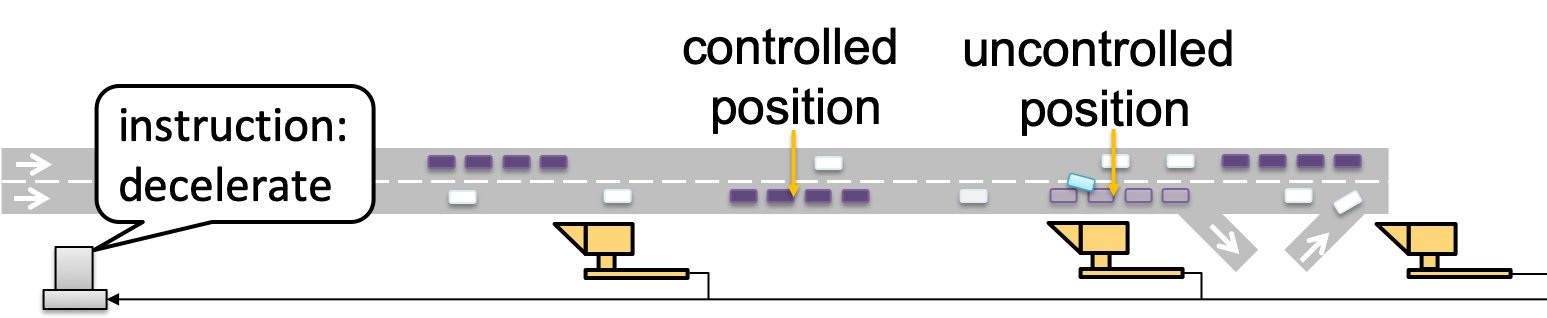}
\label{cartoon_speed}
}
\subfigure[Platoon size management.]{
\centering
\includegraphics[width=.4\textwidth]{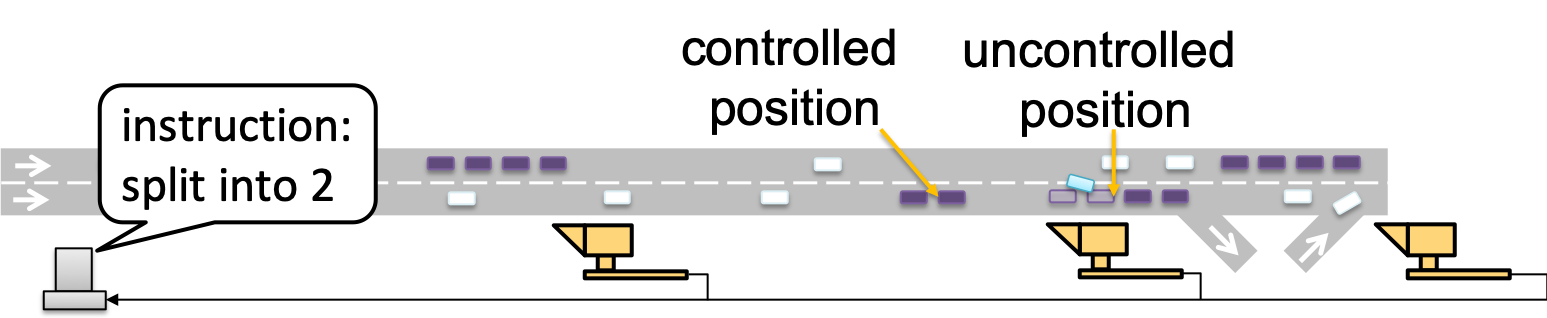}
\label{cartoon_split}
}
\caption{Two practical platoon coordination strategies.}
\label{cartoon}
\end{figure}

We now discuss how the two platoon coordination strategies presented in Section~\ref{sec_control}, viz. headway regulation and size management, can be translated to implementable instructions for platoons.
Fig.~\ref{cartoon} illustrates the implementation.
These strategies is enabled by modern vehicle-to-infrastructure (V2I) communications technologies \cite{papadimitratos2009vehicular}.
We do not explicitly consider lower-level control actions such as longitudinal and lateral control; instead, we assume that the platoons are equipped with adequate lower-level controllers that can implement the instructions from the operator.


\subsubsection{Headway regulation}
To regulate headways, the operator sends a recommended speed to each platoon when it enters the highway, and no more instructions need to be sent to this platoon; see Fig.~\ref{cartoon_speed}. The decision variable is the average speed for each platoon over the highway section, or, equivalently, the time at which a platoon is scheduled to arrive at the bottleneck.
In practice, let $L_1$ and $L_2$ be the lengths of links 1 and 2, respectively, as in Fig.~\ref{capacities}, and let $v_0$ be the nominal speed for the highway section. 
The recommended speed for an incoming platoon is
\begin{align}
    v^{hr}=\frac{L_1+L_2}{\frac{L_1+L_2}{v_0}+W^{hr}},
    \label{eq_vhr}
\end{align}
where $W^{hr}$ is given by \eqref{eq_Whr}.

\subsubsection{Size management}
To manage the size of platoons, when a platoon enters the highway, the coordination strategy first predicts the traffic condition if the platoon arrives at the bottleneck without any intervention; then, if congestion is predicted at the bottleneck, the strategy will split the platoon into shorter platoons; see Fig.~\ref{cartoon_split}.
The fluid model can be used to make such predictions.
In practice, suppose a platoon enters the highway at time $t$ and let $T$ be the solution to \eqref{eq_T}.
Specifically, the coordination decision is made in two steps:
\begin{enumerate}
\item Whether to split: The platoon is instructed to split if $Q_2^m(T)\ge\Theta-\frac{l}{2\gamma}$ and not to split otherwise; if the platoon is splitting, then the separation (i.e. headway) between the short platoons will be $l/(2\gamma(F-R-(1-\eta)\rho a))$.
\item When to arrive at bottleneck: If the platoon is not splitting, then it will travel at the nominal speed and no intervention will be needed. If the platoons are splitting into two short platoons, the leading short platoon will travel at the speed
$$
v_1^{sm}=\frac{L_1+L_2}{\frac{L_1+L_2}{v_0}+W_1^{sm}}
$$
and the following short platoon will travel at the speed
$$
v_2^{sm}=\frac{L_1+L_2}{\frac{L_1+L_2}{v_0}+W_2^{sm}}
$$
where $W^{sm}_1$ and $W^{sm}_2$ are given by \eqref{eq_W1}.
\end{enumerate}

\subsection{Simulation-based validation}
\label{sub_simulation}

The purpose of the simulations is to show that the optimal headway regulation strategy designed using the fluid model-based approach is consistent with the simulation-optimal values.
We use two standard simulation models, viz. the cell transmission model (CTM \cite{daganzo94}) and the Simulation for Urban Mobility (SUMO \cite{krajzewicz2002sumo}).
The CTM is a macroscopic traffic flow model, which evolves according to (i) the conservation law and (ii) the flow-density relation (also called ``fundamental diagram'' by transportation researchers \cite{daganzo94}).
The CTM accounts for the spatial distribution of traffic and the detailed flow-density relation, which are not captured by the fluid model.
The SUMO is a microscopic simulation model, which evolves according to vehicle-following and lane-changing behavior models for individual drivers.
Such microscopic details are captured by neither the CTM nor the fluid model.

\begin{figure}[hbt!]
\centering
\includegraphics[width=0.4\textwidth]{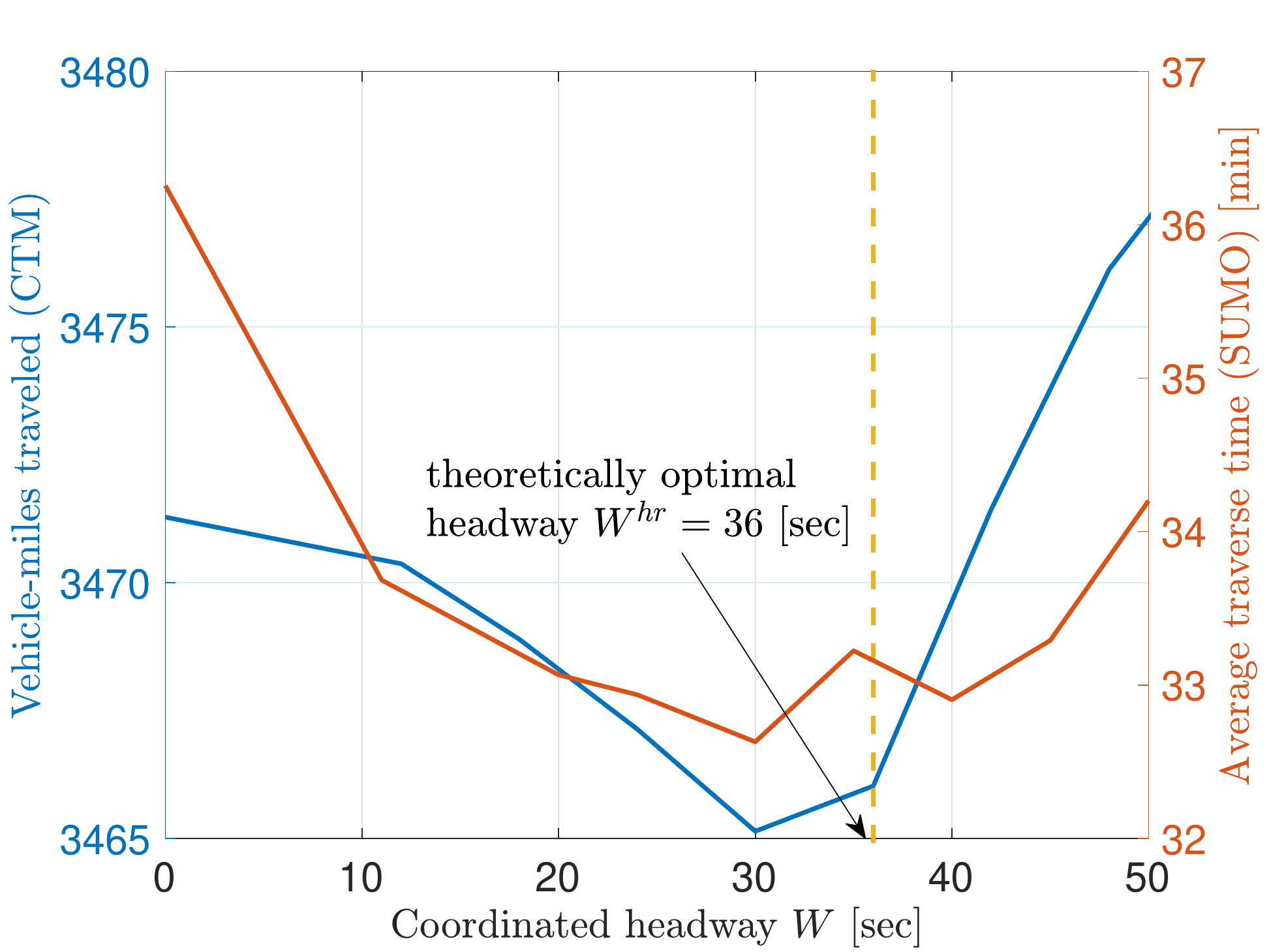}
\caption{Simulated performance metrics for various values of minimal inter-platoon headways.}
\label{fig_compare}
\end{figure}

Fig.~\ref{fig_compare} shows the simulation results with parameters in Table~\ref{tab_parameters}. The mainline demand is 2500 veh/hr, and the off-ramp demand is 1400 veh/hr.
The theoretical optimal headway $W^{hr}$ (36 sec) is close to the simulation-obtained value (30 sec).
As expected, the theoretical optimum $W^{hr}$ is greater than the simulated one. The main reason is that in the fluid model decelerating a platoon only affects the platoon itself and does not directly impact the neighboring traffic.
In both the CTM and the SUMO, however, a decelerated platoon will induce local congestion. Hence, both CTM and SUMO prefer less deceleration than the fluid model.

Next, we provide more details and discussion of the simulations.

\subsubsection{Macroscopic simulation (CTM)}
\label{sub_ctm}

We consider the cell transmission model (CTM, \cite{daganzo94}) for the highway section in Fig.~\ref{capacities}.
In particular, we consider CAVs and non-CAVs as multiple traffic classes in the CTM.
More details of the multi-class CTM is available in \cite{cicic2019coordinating}.

The parameters of the macroscopic model are chosen in accordance with the ones given in Table~\ref{tab_parameters}, i.e. the capacity of the bottleneck will be set to $F-R$, and platoon and buffer lengths chosen appropriately. Note that in order for platoons to be properly represented in this framework, we need the spatial and temporal discretization steps to be fairly short, with the platoon length spanning at least two cells. In this work, the physical platoon length is taken to be 0.1 miles, or 10 cells.

\begin{figure}[hbt!]
\centering
\subfigure[Without coordination.]{
\centering
\includegraphics[width=0.4\textwidth]{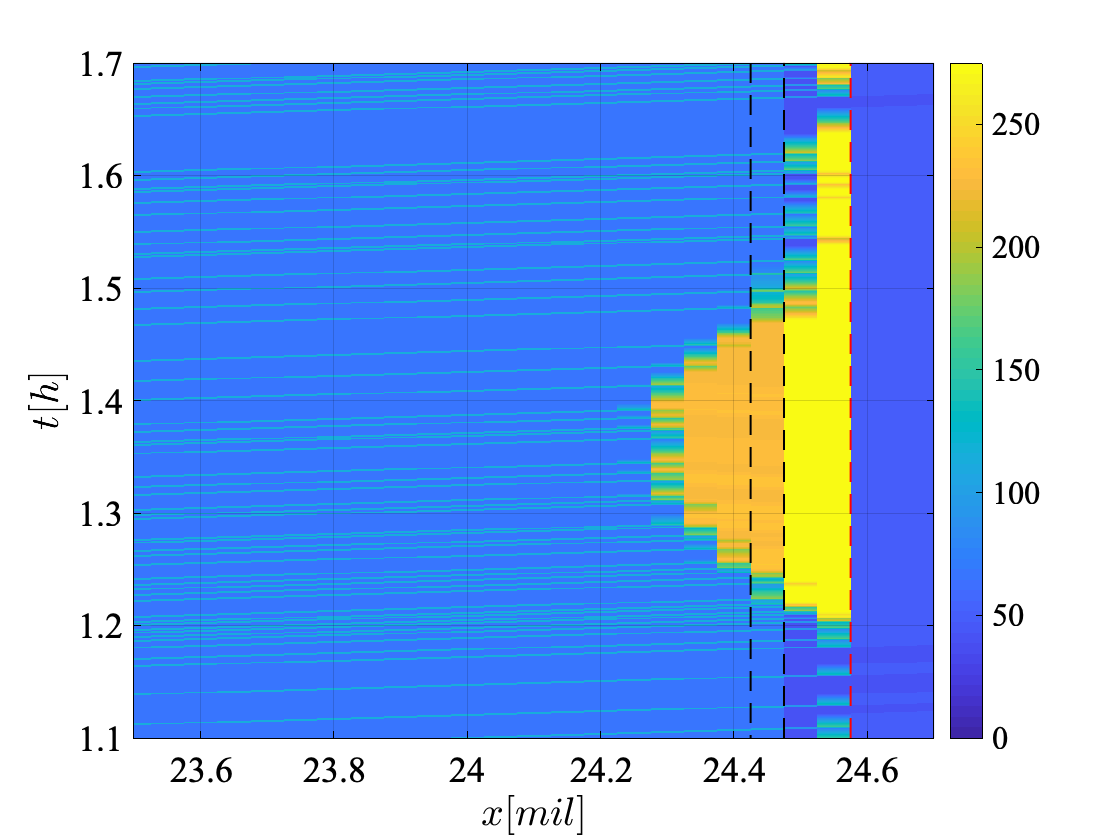}
\label{fig_sim_noctrl}
}
\subfigure[With coordination.]{
\centering
\includegraphics[width=0.4\textwidth]{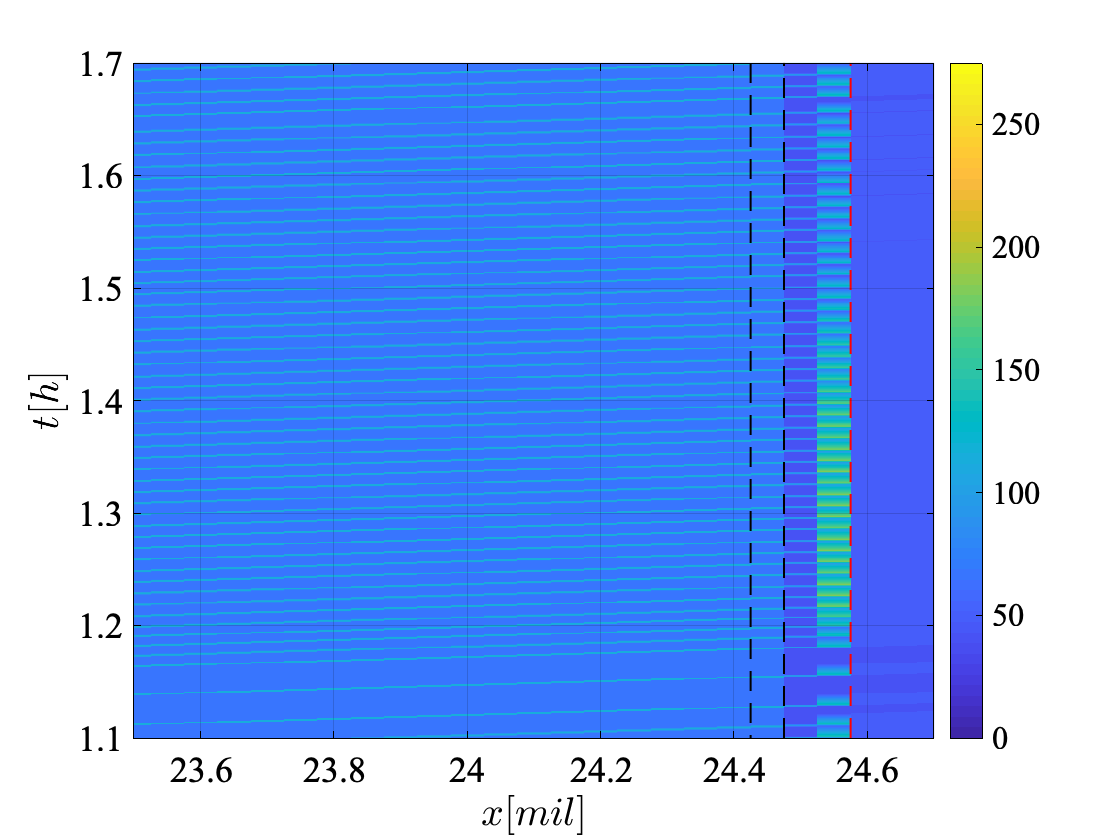}
\label{fig_sim_ctrl}
}
\caption{Traffic density contour plots for CTM simulation. The vertical red (resp. black) dashed line indicates the location of the bottleneck (resp. off-ramp). Color indicates traffic density in veh/mi.}
\label{fig_exm2}
\end{figure}

Fig.~\ref{fig_sim_noctrl} shows the simulated traffic evolution without inter-platoon coordination using color-coded traffic density.
The location of the bottleneck is shown in dashed red line, and the location of the off-ramp is outlined in dashed black line.
The streaks of brighter color represent the increased traffic density near the moving platoons.
The congestion from the bottleneck propagates upstream and the off-ramp cell becomes congested approximately one hour into the simulation.
Because of the congestion, the off-ramp is partially blocked, preventing vehicles from exiting the highway and further degrading the traffic situation.
Although the total demand $a$ is lower than the capacity of the bottleneck, randomness of the platoon arrivals may still create local congestion that disrupt the traffic flow.
Such disruptions produce congestion at the bottleneck and block the off-ramp (i.e. the bright-color area in the interval indicated by black dashed lines near distance $x=24.4$ miles), as shown in Fig.~\ref{fig_sim_noctrl} from $t \approx 1.25$ hours to $t \approx 1.5$ hours.


We apply the recommended speed \eqref{eq_vhr} to coordinate platoons in the CTM simulation.
Fig.~\ref{fig_sim_ctrl} shows the traffic evolution under headway regulation of platoons.
The simulation run considers the same situation as in the uncoordinated case (Fig.~\ref{fig_sim_noctrl}).
Whereas in the uncoordinated case the congestion from the bottleneck blocked the off-ramp, in the coordinated case we are able to spread the arrival of platoons more evenly, thus avoiding causing spillback.
Fig.~\ref{fig_compare} further illustrates how close the theoretical optimal strategy (36 sec) is closed to the simulated one (30 sec).

\begin{figure*}[hbt!]
\centering
\includegraphics[width=0.9\textwidth]{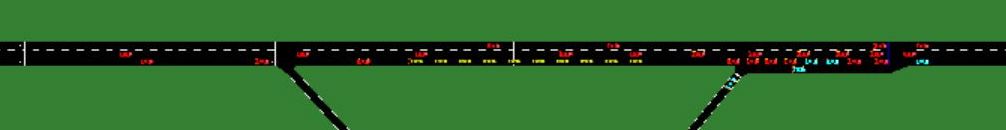}
\caption{Micro-simulation environment. Red, cyan, and yellow vehicles represent mainline traffic. on-ramp/off-ramp traffic, and CAV platoons, respectively.}
\label{fig_sumo}
\end{figure*}

\subsubsection{Microscopic simulation (SUMO)}
We also implement the headway regulation strategy introduced in Section~\ref{sub_implementation} in a micro-simulation model; see Fig.~\ref{fig_sumo}.
We use the TraCI (Traffic Control Interface) to customize the simulation and realize the functions required for this specific experiment. The TraCI features 13 individual modules varying from simulation, vehicle type, vehicle. We use Python to code the route and runner files. Some variables that we control in particular include platoon speed, total simulation time, platooning ratio, platoon length, and platooning state. We also customized the lane-changing function to prevent platoons from breaking apart at the bottleneck. The coordination instructions are realized by the runner script.

The simulation results are shown in Fig.~\ref{fig_compare}. The simulations lead to a simulation-optimal inter-platoon headway (30 sec), which is close to the value $W^{hr}$ given by the fluid model and computed via \eqref{eq_Whr} (36 sec).
A prominent pattern to be noted in Fig.~\ref{fig_compare} is that as the coordinated headway increases, the travel time in SUMO grows more slowly than the VHT in the CTM.
A reasonable explanation is that the CTM that we simulated assumes the first-come-first-serve principle: if a platoon is decelerated, then all the traffic behind will be simultaneously decelerated.
However, overtaking is allowed in SUMO, which is implemented according to SUMO's internal over-taking algorithm.
Consequently, the impact of decelerating platoons is less significant in SUMO than in the CTM.

\subsection{Further discussion}

The stochastic fluid model that we consider focuses on (i) the capacity-sharing between platoons and background traffic and (ii) the throughput gain due to reduced inter-vehicle spacing in platoons, and (iii) the impact of congestion propagation (spillback).
From a practical perspective, our model is based on the following simplifications.
\begin{enumerate}
\item The interaction between platoons and background traffic only occurs at the link boundaries. This is actually a characteristic of any queuing model. Consequently, our model does not account for interactions occurring over a distance, such as the impact of speed difference between platoons and background traffic.

\item At the interface between links 1 and 2, mainline traffic (i.e. demand $a$ and $b$) are prioritized for discharging. This is reflected in the definition of the off-ramp flow \eqref{eq_r}. Alternative models for discharging exist \cite{wright2017node}, which can be incorporated in our model as well.

\item Each link's capacity is independent of speed difference between CAVs and non-CAVs. In practice, headway regulation modifies CAVs' speeds. Consequently, highway capacity is dependent of the speed difference between CAVs and non-CAVs as well as the traffic mixture, i.e. the percentage of CAVs. In our model, the impact of heterogeneous speed can be modeled as state-dependent link capacities.

\end{enumerate}
Note that the CTM and SUMO simulations that we conducted are not restricted by the above simplifications.
The simulation results indicate that the fluid-model approach is adequate in spite of the above simplifications.
\begin{table}[hbt]
\centering
\caption{Improvement attained by various coordination strategies.}
\label{tab_simulation}
\begin{tabular}{|l|c|c|c|c|}
\hline
Strategy                                                        & \begin{tabular}[c]{@{}l@{}}VHT,\\ CTM\end{tabular} & \begin{tabular}[c]{@{}l@{}}Improve-\\ ment (VHT)\end{tabular} & \begin{tabular}[c]{@{}l@{}}Traverse time,\\ SUMO {[}min{]}\end{tabular} & \begin{tabular}[c]{@{}l@{}}Improve-\\ ment {[}min{]}\end{tabular} \\ \hline
\begin{tabular}[c]{@{}l@{}}No\\ coordination\end{tabular}       & 3471                                               & 0                                                             & 36.26                                                                   & 0                                                                 \\ \hline
\begin{tabular}[c]{@{}l@{}}Theoretically\\ optimal\end{tabular} & 3466                                               & 5.3                                                           & 33.27                                                                   & 2.99                                                              \\ \hline
\begin{tabular}[c]{@{}l@{}}Simulation\\ optimal\end{tabular}   & 3465                                               & 6.2                                                           & 32.58                                                                   & 3.68                                                               \\ \hline
\end{tabular}
\end{table}
As Table~\ref{tab_simulation} shows, compared to the baseline scenario without headway regulation, the theoretically optimal strategy attains 85\% (resp. 81\%) of the improvement attained by the simulation optimum in CTM (resp. SUMO).

%% file: Texts/5_Conclusion.tex
\section{Concluding remarks}
\label{sec_conclude}

In this paper, we develop a stochastic fluid model for analysis and design of platoon coordination strategies. The model focuses on the interaction between CAV platoons and non-CAVs, the impact of key platooning parameters (including platooning ratio and platoon size), and road geometry (buffer space). Based on the theory of Markov processes and queuing properties of the model, we derive theoretical bounds for the throughput of the model and identify a set of coordination strategies that maximize throughput as well as minimize delay. We discuss how such strategies can be implemented in practice and validate the fluid model-based results 
using standard macroscopic and microscopic simulation environments.
Our results are useful for link-level coordination of CAV platoons.

This work can be extended in the following directions: (i) evaluation of macroscopic impact due to various vehicle-level controllers, (ii) extension to multiple-section highways, and (iii) integration with network-level scheduling and routing of CAVs.

\section*{Appendix 1. Proof of Theorem~\ref{thm_sufficient}}
\label{sub_sufficient}

The stability criterion is obtained by showing that the Lyapunov function
\begin{align}
    V(q):=&\frac12(q_0^m+q_1^m+q_2^m)^2\nonumber\\
    &+\bigg(k\Big(q_0^m+q_1^m+\frac1{2\Theta}(q_2^m)^2\Big)+q_1^o\bigg)q_1^o
    \label{eq_V}
\end{align}
satisfies the Foster-Lyapunov criterion if \eqref{eq_a<a*}--\eqref{eq_max} hold.
Specifically, we show that there exist $k>0$, $c>0$, and $d<\infty$ such that
\begin{align}
    \mathscr LV(q)\le-c|q|+d,
    \quad\forall q\in{\mathcal M},
    \label{eq_LV}
\end{align}
which implies stability via the Foster-Lyapunov criterion.
To proceed, we decompose $V$ into
\begin{subequations}\begin{align}
    &V^m(q):=\frac12(q_0^m+q_1^m+q_2^m)^2,\label{eq_Vm}\\
    &V^o(q):=\bigg(k\Big(q_1^m+\frac1{2\Theta}(q_2^m)^2\Big)+q_1^o\bigg)q_1^o\label{eq_Vo}
\end{align}\end{subequations}
and show that there exists $c^m,c^o>0$ and $d^m,d^o<\infty$ such that for all $q\in\mathcal M_{\mu,\nu}$,
\begin{subequations}
\begin{align}
    &\mathscr LV^m(q)\le-c^m(q_0^m+q_1^m+q_2^m)+d^m,\label{eq_LVa}\\
    &\mathscr LV^o(q)\le-c^oq_1^o+d^o.\label{eq_LVb}
\end{align}
\end{subequations}
Note that the above implies \eqref{eq_LV} and hence stability.
The rest of this subsection is devoted to the proof of \eqref{eq_LVa}--\eqref{eq_LVb}.

\subsubsection{Proof of \eqref{eq_LVa}}
For $q$ such that $q_2^m=0$, we have
\begin{align*}
    \mathscr LV^m(q)=\frac12(\lambda l/\gamma)^2=\frac{\eta \rho a}{2\gamma}
\end{align*}
For $q$ such that $q_2^m>0$, we have
\begin{align*}
    \mathscr LV^m(q)&=((1-\eta)\rho a+\lambda l/\gamma-(F-R))(q_0^m+q_1^m+q_2^m)\\
    &=((\eta/\gamma+1-\eta)\rho a-(F-R))(q_0^m+_1^m+q_2^m).
\end{align*}
Hence, we have 
\begin{align*}
    &c^m=(F-R)-(\eta/\gamma+1-\eta)\rho a
    \stackrel{\footnotesize\eqref{eq_a<a*}}>0,\\
    &d^m=\frac{\eta \rho a}{2\gamma}<\infty
\end{align*}
that satisfy \eqref{eq_LVa}.

\subsubsection{Proof of \eqref{eq_LVb}}
\label{subsub_LVo}

For $q$ such that $q_1^o=0$, we have
\begin{align}
    \mathscr LV^o(q)=0.
    \label{eq_LVb0}
\end{align}
For $q$ such that $q_1^o>0$, we need to consider the following cases:
\begin{enumerate}
    \item $q_0^m=q_1^m=0$. In this case, we require a $c^o>0$ such that
    \begin{align}
        &\mathscr LV^o(q)\nonumber\\
        &= \bigg(k(q_2^m/\Theta)((1-\eta)\rho a-(F-R))+\lambda k\Big(S^m_0(q;\nu(q))\nonumber\\
        &\quad+S^m_1(q;\nu(q))+\frac{1}{2\Theta}(S^m_2(q;\nu(q)))^2-\frac1{2\Theta}(q_2^m)^2\Big)\nonumber\\
        &\quad+(1-\rho)a-R\bigg)q_1^o\nonumber\\
        &\le-c^oq_1^o,
        \quad \forall q:q_1^m=0,q_1^o>0.
        \label{eq_LVb1}
    \end{align}
    By Assumption~\ref{asm_nu}, we have
    \begin{align*}
        &S^m_0(q;\nu(q))+S^m_1(q;\nu(q))+\frac{1}{2\Theta}(S^m_2(q;\nu(q)))^2\\
        &-\frac1{2\Theta}(q_2^m)^2\le S^m_0(q;\nu(q))+S^m_1(q;\nu(q))\\
        &+\frac{1}{2\Theta}(S^m_2(q;\nu(q)))^2-\frac1{2\Theta}(q_2^m)^2\Big|_{q_1^o=0},
        \quad\forall q\in\mathcal Q.
    \end{align*}
    Hence, we require a $k>0$ such that
    \begin{align}
        &(q_2^m/\Theta)((1-\eta)\rho a-(F-R))+\lambda \Big(S^m_0(q;\nu(q))\nonumber\\
        &+S^m_1(q;\nu(q))+\frac{1}{2\Theta}(S^m_2(q;\nu(q)))^2-\frac1{2\Theta}(q_2^m)^2\Big|_{q_1^o=0}\Big)\nonumber\\
        &<\frac{R-(1-\rho)a}k,
        \quad \forall q:q_0^m=q_1^m=0,q_1^o>0.
        \label{eq_Db1}
    \end{align}
    
    \item $q_0^m>0$ or $q_1^m>0,\ q_2^m<\Theta$. In this case, we require $c^0>0$ such that
    \begin{align}
        &\mathscr LV^o(q)=\bigg(k\Big((1-\eta)\rho a-F+(q_2^m/\Theta)(F\nonumber\\
        &\quad-(F-R))\Big)+k\lambda(l/\gamma)\Big)+(1-\rho)a-R\bigg)q_1^o\nonumber\\
        &=\Big(k\Big((\frac{\eta}{\gamma}+1-\eta)\rho a-(F-R)\Big)+(1-\rho)a-R\Big)q_1^o\nonumber\\
        &\le-c^oq_1^o,
        \quad \forall q:q_0^m>0\mbox{ or }q_1^m>0,q_1^o<\Theta.
        \label{eq_LVb2}
    \end{align}
    For all $k>0$, the existence of such a $c^o$ is ensured by \eqref{eq_a<a*}.
    
    \item $q_0^m>0$ or $q_1^m>0,\ q_2^m=\Theta$. In this case,
    \begin{align}
        \mathscr LV^o(q)&=\bigg(k\Big((\eta/\gamma+1-\eta)\rho a-(F-R)\Big)\nonumber\\
        &\quad+(1-\rho)a\bigg)q_1^o\le-c^oq_1^o.
        \label{eq_LVb3}
    \end{align}
    We require a $k>0$ such that
    \begin{align}
       k\Big((F-R)-(\eta/\gamma+1-\eta)\rho a\Big)>(1-\rho)a.
       \label{eq_Db3}
    \end{align}
\end{enumerate}
Note that \eqref{eq_max} ensures the existence of $k>0$ that simultaneously satisfies \eqref{eq_Db1} and \eqref{eq_Db3}. Hence, there exists $k>0$ and $c^o>0$ such that satisfy \eqref{eq_LVb1}, \eqref{eq_LVb2}, and \eqref{eq_LVb3}, which, together with \eqref{eq_LVb0}, imply \eqref{eq_LVb}.

\section*{Appendix 2. Proof of Theorem~\ref{thm_bounds}}
\label{sub_bounds}



\emph{1) Lower bound}
The lower bound is the minimum of the following two terms:
\begin{align*}
    &\underline a_1:=\frac{F-R}{\rho(\eta/\gamma+1-\eta)},\\
    &\underline a_2:=
\frac{R}{1-\rho+\frac{1}{2}\Big(\sqrt{\zeta^2+\frac{2\rho Rl}{\gamma\Theta(F-R)}}-\zeta)\Big)},
\end{align*}
where $\zeta$ is given by \eqref{eq_zeta0}.
We prove the lower bound by applying Theorem~\ref{thm_sufficient}, i.e. verifying that \eqref{eq_a<a*}--\eqref{eq_max} hold if $a<\min\{\underline a_1,\underline a_2\}$:
\begin{enumerate}
    \item Since $2\rho Rl/(\gamma\Theta(F-R))>0$, we have $\underline a_2\le R/(1-\rho)$. Hence, $a<\underline a_2$ and $a<\underline a_1$ ensure that \eqref{eq_a<a*} holds.
    \item For $\xi\in\mathcal Q$ such that $\xi_0^m=\xi^m_1=0$,
\begin{align*}
    &\frac{\xi_2^m}{\Theta}\Big((1-\eta)\rho a-(F-R)\Big)+\frac{\eta\rho a}{l} \bigg(S^m_0(\xi;0)\nonumber\\
    &\qquad+S^m_1(\xi;0)+\frac{1}{2\Theta}\Big(S^m_2(\xi;0)\Big)^2-\frac1{2\Theta}(\xi_2^m)^2\bigg)\\
    &=\Big(\frac{\xi_2^m}{\Theta}((1-\eta)\rho a-(F-R)\Big)\\
    &\qquad+\frac{\lambda }{2\Theta}(2\xi_2^m l/\gamma+(l/\gamma)^2).
\end{align*}
Since $a<\underline a_1$ and since $\xi_2^m\le\Theta$, the above implies that
\begin{align*}
    &\max_{\substack{\xi\in\mathcal Q:\xi_0^m\\=\xi_1^m=\xi_1^o=0}}\frac{\xi_2^m}{\Theta}\Big((1-\eta)\rho a-(F-R)\Big)+\frac{\eta\rho a}{l} \bigg(S^m_0(\xi;0)\nonumber\\
    &\qquad+S^m_1(\xi;0)+\frac{1}{2\Theta}\Big(S^m_2(\xi;0)\Big)^2-\frac1{2\Theta}(\xi_2^m)^2\bigg)\\
    &\le \Big(\frac{\lambda k}{2\Theta}(2\Theta l/\gamma+(l/\gamma)^2).
\end{align*}
With the above, one can verify that if $a<\underline a_2$, then \eqref{eq_max} holds.
\end{enumerate}

\emph{2) Upper bound}

We prove the upper bound by showing that if $\{Q(t);t>0\}$ is stable, then $a\le\min\{\bar a_1,\bar a_2\}$, where
\begin{align*}
    &\bar a_1:=\frac{F-R}{\rho(\eta/\gamma+1-\eta)},\\
    &\bar a_2:=\frac{(1-\omega)R}{1-\rho}.
\end{align*}

\begin{enumerate}
    \item $a\le\bar a_1$ can be obtained from the nominal throughput given by Lemma~\ref{lmm_maximum}.
    \item To show $a\le a_2$, note that when $a<\bar a_1$, the M/D/1 process $\{N(t);t>0\}$ is stable and admits a steady-state distribution $\{\pi_n;n=0,1,\ldots\}$ defined in \eqref{eq_pin}.
Hence, there exists $\omega_0$ and $\omega$ such that
\begin{align*}
    &\lim_{t\to\infty}\frac1t\int_{\tau=0}^t\mathbb I_{Q_1^m>0,Q_2^m=\Theta}d\tau=\omega_0\quad a.s.,\\
    &\lim_{t\to\infty}\frac1t\int_{\tau=0}^t\mathbb I_{N(t)\ge\ceil{\gamma\Theta/l}}d\tau=\omega\quad a.s.
\end{align*}
where $\omega$ is in fact given by \eqref{eq_omega}.

Next, consider the set $\mathcal M_{0,0}\subset\mathcal Q$ defined by
\begin{align}\label{eq_M}
{\mathcal M_{0,0}}=\Big((\{0\}\times[0,\Theta])\cup((0,\infty)\times\{\Theta\})\Big)\times[0,\infty).
\end{align}
One can show that $\mathcal M_{0,0}$ is an invariant set.
Hence, for each initial condition $q\in\mathcal M_{\mu,\nu}$, we have $Q(t)\in\mathcal M_{\mu,\nu}$ for all $t>0$.
Thus we have $Q_2^m(t)=\Theta$ if $Q_1^m(t)>0$ for sufficiently large $t$. 
Hence, if $Q_1^m(t)+Q_2^m(t)>\Theta$, i.e. if $Q_1^m(t)>0$ and $Q_2^m(t)=\Theta$, then $N(t)\ge\ceil{\gamma\Theta/l}$.
Therefore, we have
$
\omega_0\ge\omega.
$
Finally, note that if $|Q(t)|$ is bounded, then 
\begin{align*}
    (1-\rho)a&\le\lim_{t\to\infty}\frac1t\int_{\tau=0}^tr(Q(\tau);0)d\tau\\
    &=\lim_{t\to\infty}\frac1t\Big(\int_{\tau:Q_1^m(\tau)+Q_1^m(\tau)\le\Theta}^tr(Q(\tau);0)d\tau\\
    &\quad+\int_{\tau:Q_1^m(\tau)+Q_1^m(\tau)>\Theta}^tr(Q(\tau);0)d\tau\Big)\\
    &\stackrel{a.s.}=(1-\omega_0)R
    \le(1-\omega)R=\bar a_2.
\end{align*}
\end{enumerate}

\section*{Appendix 3. Proof of Theorem~\ref{thm_optimal}}
\label{sub_optimal}

\emph{1) Stability}
\label{sub_max}

The necessity of \eqref{eq_iff} results from Lemma~\ref{lmm_maximum}.
To show the sufficiency, one can indeed use Theorem~\ref{thm_sufficient} to show that the model is stabilized by $(\mu,\nu)\in\mathscr U^*\times\mathscr V^*$ if $a<a^*$ in the sense of a bounded 1-norm.
In this subsection, we use an alternative Lyapunov function
\begin{align}
\tilde V(q)=e^{\beta|q|},
\quad q\in\mathcal Q
\label{eq_Vtilde}
\end{align}
and obtain a stronger stability:
\begin{align}
    \limsup_{t\to\infty}\frac1t\int_{\tau=0}^t{\mathrm E}[e^{\beta|{Q}(\tau)|}]d\tau\le Z.
    \label{eq_MGF}
\end{align}
That is, the state is bounded in the MGF.

To proceed, define
\begin{align}
\zeta_{\mu,\nu}:=&\max\Big\{\max_{q\in\mathcal Q}S_2^m(q;\nu(q)),\ \sup\{\zeta\ge0:\nonumber\\
&(\forall {q}:q_2^m=\zeta)\ {\mu}(q)\ge {F-R}-(1-\eta)\rho a\}\Big\}.
\label{eq_zeta}
\end{align}
We know from \eqref{eq_muq=0}--\eqref{eq_S1m} that $\zeta_{\mu,\nu}<\Theta$ for all $(\mu,\nu)\in\mathscr U^*\times \mathscr V^*$.
Then, consider the set
\begin{align}\label{eq_Mtilde}
\mathcal M_{\mu,\nu}=[0,\infty)\times\{0\}^2\times[0,\zeta_{\mu,\nu}].
\end{align}

\begin{lmm}
The set $\mathcal M_{\mu,\nu}$ as defined in \eqref{eq_Mtilde} is an invariant set.
\end{lmm}
\noindent\emph{Proof.}
%
(i) For all $q\in\mathcal M_{\mu,\nu}$, we have
\begin{align*}
    G^m_1(q;\mu(q))&=\mu(q)+((1-\eta)\rho+(1-\rho))a\\
    &\quad-\min\Big\{\Big((1-\eta)\rho+(1-\rho)\Big)a+\mu(q),F\Big\}\\
    &\stackrel{\footnotesize\eqref{eq_muq<F}}=0.
\end{align*}
Hence, $\{0\}$ is invariant for $q_1^m$.
In addition, for $q\in\mathcal M_{\mu,\nu}$ such that $q_2^m=\zeta$, we have
\begin{align*}
    G_2^m(q;\mu(q))&=f_1(q;\mu(q))-f_2(q;\mu(q))\\
    &\stackrel{\footnotesize\eqref{eq_zeta}}<(\mu(q)+(1-\eta)\rho a)-(\mu(q)+(1-\eta)\rho a)\\
    &<0.
\end{align*}
Hence, $[0,\zeta_{\mu,\nu}]$ is invariant for $q_2^m$.

(ii) By \eqref{eq_S1m}, since $S_1^m(q;\nu(q))=0$ for all $q\in\mathcal Q$, $\{0\}$ is invariant for $q_1^m$. By \eqref{eq_zeta}, $S_2^m(q;\nu(q))\le\zeta$ for all $q\in\mathcal Q$. Hence, $[0,\zeta]$ is invariant for $q_2^m$.
\hfill$\square$

In the rest of this proof, we show that there exists $\beta>0$, $c>0$, and $\tilde d<\infty$ verifying the drift condition
\begin{align}\label{eq_drift2}
{\mathcal L}\tilde V(q)\le-ce^{\beta|q|}+ \tilde d
\quad\forall q\in\mathcal M_{\mu,\nu},
\end{align}
which leads to \eqref{eq_MGF} by the Foster-Lyapunov criterion.

We partition the invariant set $\mathcal M_{\mu,\nu}$ into two subsets:
\begin{align*}
&{\mathcal M_{\mu,\nu}^0}=\{0\}^3\times[0,\zeta_{\mu,\nu}],\\
&{\mathcal M_{\mu,\nu}^1}=(0,\infty)\times\{0\}^2\times[0,\zeta_{\mu,\nu}],
\end{align*}
For $q\in{\mathcal M_{\mu,\nu}^0}$, we have
\begin{align*}
{\mathcal L}\tilde V(q)
&=\nabla_{{q}}e^{\beta|q|}G(q;\mu(q))+\lambda(e^{\beta |S(q;\nu(q))|}-e^{\beta|q|})\\
&=(\beta((1-\eta)\rho a-\tilde f_2(q;\mu(q)))+\lambda(e^{\beta l/\gamma}-1))e^{\beta|q|}\\
&\le\lambda(e^{\beta l/\gamma}-1)e^{\beta l/\gamma}=\tilde d^*,
\end{align*}
which also defines $\tilde d^*$.
For $q\in{\mathcal M_{\mu,\nu}^1}$, we have
\begin{align*}
{\mathcal L}\tilde V(q)
&
=(\beta((1-\eta)\rho a-(F-R))+\lambda(e^{\beta l/\gamma}-1))e^{\beta|q|}\\
&=\phi(\beta)e^{\beta|q|},
\end{align*}
where the definition of the function $\phi$ is clear.
Since
\begin{align*}
\phi(0)=0,\
\frac{d}{d\beta}\phi(\beta)\Big|_{\beta=0}=(1-\eta)\rho a-(F-R)+\eta\rho a/\gamma,
\end{align*}
there exists $\beta^*>0$ such that $\phi(\beta^*)<0$ if $a<a^*$.

In conclusion, there exist $\beta=\beta^*,$ $c=\phi(\beta^*)$, and $\tilde d=\tilde d^*$ that verify \eqref{eq_drift2}. Then, by the Foster-Lyapunov stability criterion, we conclude \eqref{eq_MGF}.


\emph{2) Queue minimization}

Next, we use a sample path-based method to show that any $(\mu,\nu)\in\mathscr U^*\times\mathscr V^*$ minimizes the total queue size $|Q(t)|$ for all $t$ over all control policies $(\mu,\nu)\in\mathscr U\times\mathscr V$.

Let $\{{M}(t);t>0\}$ be the counting process of platoon arrivals.
For a given sample path $\{{m}(t);t>0\}$ of the counting process and a given initial condition $q\in{\mathcal Q}$, let $\{{q}(t);t>0\}$ and $\{\psi(t);t>0\}$ be the corresponding trajectories under a control policy $(\mu^*,\nu^*)\in\mathscr U^*\times\mathscr V^*$ and under a control policy $(\mu,\nu)\in\mathscr U\times\mathscr V$, respectively.
To show the optimality of $(\mu^*,\nu^*)$, it suffices to show $|{q}(t)|\le|\psi(t)|$ for any $t\ge0$.
Without loss of generality, we consider zero initial condition. 
We prove this by contradiction as follows.

Assume by contradiction that there exists $(\mu,\nu)\in\mathscr U\times\mathscr V$ such that
\begin{align}
\exists t_1>0,\ |q(t_1)|>|\psi(t_1)|.\label{eq_contradiction}
\end{align}
Between resets, the continuity of ${q}(t)$ and $\psi(t)$ follows from Assumption~\ref{asm_mu}.
Therefore, there must exist a ``crossing time'' $t_0\in[0,t)$ such that
\begin{align}
|{q}(t_0)|=|\psi(t_0)|,\
\frac d{dt}|q(t)|>\frac d{dt}|\psi(t)|.
\label{crossing}
\end{align}
Note that the ``crossing'' must happen between resets. To see this, recall that Assumption~\ref{asm_nu} ensures that if $|{q}(t_-)|=|\psi(t_-)|$ and if a reset occurs at time $t$, then
\begin{align*}
    |q(t)|=|\psi(t)|=|q(t_-)|+l/\gamma.
\end{align*}

Since the system admits the invariant set $\mathcal M_{\mu^*,\nu^*}$ as given in \eqref{eq_Mtilde} under $(\mu^*,\nu^*)$, we have ${q}_1^m(t)={q}_1^o(t)=0$ for all $t\ge0$.
Hence, a necessary condition for \eqref{crossing} is that
\begin{subequations}
\begin{align}
&{q}_0^m(t_0)+{q}_2^m(t_0)
\ge\psi_0^m(t_0)+\psi_1^m(t_0)+\psi_2^m(t_0),\\
&G_0^m(q(t_0);\mu^*)+G_2^m(q(t_0);\mu^*)\nonumber\\
&>G_0^m(\psi(t_0);\mu)+G_1^m(\psi(t_0);\mu)+G_2^m(\psi(t_0);\mu)).\label{eq_dot*=dotd}
\end{align}
\end{subequations}
However, one can obtain from \eqref{eq_G0}--\eqref{eq_G2} that if ${q}_0^m(t_0)+{q}_2^m(t_0)
\ge\psi_0^m(t_0)+\psi_1^m(t_0)+\psi_2^m(t_0)=0$, then
\begin{align}
&G_0^m(q(t_0);\mu^*)+G_2^m(q(t_0);\mu^*)\nonumber\\
&
=G_0^m(\psi(t_0);\mu)+G_1^m(\psi(t_0);\mu)+G_2^m(\psi(t_0);\mu));\label{eq_sum=0}
\end{align}
if ${q}_0^m(t_0)+{q}_2^m(t_0)
\ge\psi_0^m(t_0)+\psi_1^m(t_0)+\psi_2^m(t_0)>0$, then
\begin{subequations}
\begin{align}
&G_0^m(q(t_0);\mu^*)+G_2^m(q(t_0);\mu^*)=(1-\eta)\rho a-(F-R),\label{eq_sum>0_1}\\
&G_0^m(\psi(t_0);\mu)+G_1^m(\psi(t_0);\mu)+G_2^m(\psi(t_0);\mu))\nonumber\\
&\quad=(1-\eta)\rho a-r(\psi(t_0))
\ge (1-\eta)\rho a-(F-R).\label{eq_sum>0_2}
\end{align}
\end{subequations}
Since both \eqref{eq_sum=0} and \eqref{eq_sum>0_1}--\eqref{eq_sum>0_2} contradict with \eqref{eq_dot*=dotd}, we conclude that $(\mu,\nu)$ cannot achieve \eqref{eq_contradiction}. That is, if $|q(t_1)|=|\psi(t_1)|$, then $|q(t)|$ cannot increase faster (or decrease slower) than $|\psi(t)|$ at time $t=t_1$. This proves the optimality of $(\mu^*,\nu^*)$.

\emph{3) Mean queue size}
\label{sub_queue}

To show that the 1-norm of the state converges as in \eqref{eq_Qbar}, we first need to show that the process $\{Q(t);t\ge0\}$ is ergodic. Formally, let $P_t(q)$ be the distribution of $Q(t)$ given the initial condition $Q(0)=q$ for $q\in\mathcal Q$. Then, there exists a unique probability measure $P^*$ on $\mathcal Q$ such that
\begin{align}
    \lim_{t\to\infty}\|P_t(q)-P^*\|_{\mathrm{TV}}=0,
    \label{eq_TV}
\end{align}
where $\|\cdot\|_{\mathrm{TV}}$ is the total-variation distance between two probability measures \cite{benaim15}. Ergodicity ensures convergence of the time average towards the expected value, if the expected value exists.

\emph{Convergence}:
The Foster-Lyapunov criterion ensures the existence of an invariant measure $P^*$ \cite[Theorem 4.5]{meyn93}.
We only need to further show that the invariant measure is unique.
This can be shown via the ``coupling'' condition:

\noindent{\bf Coupling condition \cite{meyn93}.}
\emph{Let $q,q'\in\mathcal Q$ be two initial conditions and $Q(t),Q'(t)$ be the trajectories starting therefrom. Then, there exists $\delta>0$ and $T<\infty$ such that
\begin{align}
    \Pr\{Q(T)=Q'(t)|Q(0)=q,Q'(0)=q'\}=\delta.
    \label{eq_coupling}
\end{align}
}

To show that the stochastic fluid model controlled by $(\mu,\nu)\in\mathscr U^*\times\mathscr V^*$ satisfies the coupling condition, note that for an arbitrary initial condition $q\in{\mathcal Q}$, there exists
$
T=(\sum_{k=1}^3q_k^m)/(F-R)+q_1^o/R
$
such that
\begin{align*}
\Pr\{Q(T)=0|Q(0)=q\}\ge e^{-\lambda T}>0.
\end{align*}
Then, by \cite[Theorem 6.1]{meyn1993stability}, the above, together with \eqref{eq_drift2}, implies convergence in the sense of \eqref{eq_TV}.

\begin{rem}
In fact, the above argument ensures {exponentially convergent} \cite{meyn93} in the sense that there exist a constant $\kappa>0$ and a finite-valued function $U:\tilde{\mathcal Q}\to\mathbb R_{\ge0}$ such that
$$
\|P_t(q)-P^*\|_{\mathrm{TV}}\le U(q)e^{-\kappa t}
\quad\forall t\ge0.
$$
\end{rem}

\begin{figure}[hbt]
\centering
\includegraphics[width=0.3\textwidth]{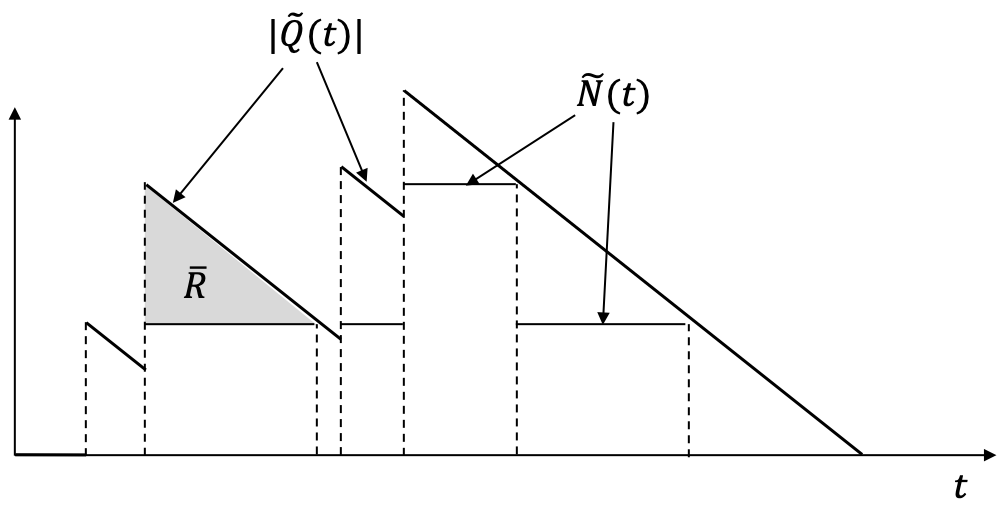}
\caption{The controlled fluid process ${Q}(t)$ envelops an M/D/1 process $\tilde N(t)$.}
\label{subprocesses}
\end{figure}

\emph{Queuing delay}:
The evolution of the total queue length $|Q(t)|$ can be viewed as the superposition of two subprocesses; see Fig.~\ref{subprocesses}. 
\begin{enumerate}
    \item The first process is an M/D/1 process $\{N(t);t\ge0\}$ defined by \eqref{eq_N}.
By the Pollazcek-Khinchin formula \cite[p. 248]{gallager13} and the Little's theorem~\cite[Theorem 5.5.9]{gallager13}, the mean number of waiting jobs (excluding the one being served) of this process is
\begin{align*}
&\bar N=\frac{\lambda^2 s^2}{2(1-\lambda s)}\\
&=\frac{\eta^2\rho^2 a^2/(2\gamma^2)}{(F-R-(1-\eta)\rho a)(F-R-(\eta/\gamma+1-\eta)\rho a)}.
\end{align*}
    \item The second process is the ``services'' each job experiences. As Fig.~\ref{subprocesses} shows, the cumulative queuing delay during services is given by
$$
\bar R=\frac{l^2}{2\gamma^2(F-R-(1-\eta)\rho a)}.
$$
Thus, the total fluid queue length is
$$
\bar Q=\frac l\gamma\bar N+\frac{\eta\rho a}{l}\bar R,
$$
which leads to \eqref{eq_Qbar}.
\end{enumerate}